\documentclass[11pt,reqno]{amsart}
\usepackage{amsthm}
\usepackage{amssymb}
\usepackage{latexsym}
\usepackage{multicol}
\usepackage{verbatim,enumerate}
\usepackage{hyperref}
\usepackage{amsmath, amscd}

\advance\textwidth by 1.2in \advance\oddsidemargin by -.6in \advance\evensidemargin by -.6in
\usepackage{tikz}
\usetikzlibrary{decorations.markings}
\tikzstyle{vertex}=[ circle, fill, draw, inner sep=0pt, minimum size=4pt,]
\tikzstyle{edge}= [thick]

\newtheorem*{cor}{Corollary}%[section]
\newtheorem*{lem}{Lemma}
\newtheorem*{prop}{Proposition}

\theoremstyle{definition}
\newtheorem*{defn}{Definition}
\theoremstyle{definition}
\newtheorem{thm}{Theorem}
\newtheorem*{thm*}{Theorem}

\newtheorem*{rem}{Remark}

\newenvironment{pf}{\proof}{\endproof}
\newcounter{cnt}
\newenvironment{enumerit}{\begin{list}{{\hfill\rm(\roman{cnt})\hfill}}{%
\settowidth{\labelwidth}{{\rm(iv)}}\leftmargin=\labelwidth%
\advance\leftmargin by \labelsep\rightmargin=0pt\usecounter{cnt}}}{\end{list}} \makeatletter
\def\mydggeometry{\makeatletter\dg@YGRID=1\dg@XGRID=20\unitlength=0.003pt\makeatother}
\makeatother \theoremstyle{remark}

\numberwithin{equation}{section}

\let\bwdg\bigwedge
\def\bigwedge{{\textstyle\bwdg}}

\begin{document}

\newcommand{\thmref}[1]{Theorem~\ref{#1}}
\newcommand{\secref}[1]{Section~\ref{#1}}
\newcommand{\lemref}[1]{Lemma~\ref{#1}}
\newcommand{\propref}[1]{Proposition~\ref{#1}}
\newcommand{\corref}[1]{Corollary~\ref{#1}}
\newcommand{\remref}[1]{Remark~\ref{#1}}
\newcommand{\defref}[1]{Definition~\ref{#1}}
\newcommand{\er}[1]{(\ref{#1})}
\newcommand{\id}{\operatorname{id}}
\newcommand{\ord}{\operatorname{\emph{ord}}}
\newcommand{\sgn}{\operatorname{sgn}}
\newcommand{\wt}{\operatorname{wt}}
\newcommand{\tensor}{\otimes}
\newcommand{\from}{\leftarrow}
\newcommand{\nc}{\newcommand}
\newcommand{\rnc}{\renewcommand}
\newcommand{\dist}{\operatorname{dist}}
\newcommand{\qbinom}[2]{\genfrac[]{0pt}0{#1}{#2}}
\nc{\cal}{\mathcal} \nc{\goth}{\mathfrak} \rnc{\bold}{\mathbf}
\renewcommand{\frak}{\mathfrak}
\newcommand{\supp}{\operatorname{supp}}
\newcommand{\Irr}{\operatorname{Irr}}
\newcommand{\psym}{\mathcal{P}^+_{K,n}}
\newcommand{\psyml}{\mathcal{P}^+_{K,\lambda}}
\newcommand{\psymt}{\mathcal{P}^+_{2,\lambda}}
\renewcommand{\Bbb}{\mathbb}
\nc\bomega{{\mbox{\boldmath $\omega$}}} \nc\bpsi{{\mbox{\boldmath $\Psi$}}}
 \nc\balpha{{\mbox{\boldmath $\alpha$}}}
 \nc\bpi{{\mbox{\boldmath $\pi$}}}
  \nc\bxi{{\mbox{\boldmath $\xi$}}}
\nc\bmu{{\mbox{\boldmath $\mu$}}} \nc\bcN{{\mbox{\boldmath $\cal{N}$}}} \nc\bcm{{\mbox{\boldmath $\cal{M}$}}} \nc\blambda{{\mbox{\boldmath
$\lambda$}}}%\nc\mathbb Nu{{\mbox{\boldmath $\nu$}}}

\newcommand{\Tmn}{\bold{T}_{\lambda^1, \lambda^2}^{\nu}}

\newcommand{\lie}[1]{\mathfrak{#1}}
\newcommand{\ol}[1]{\overline{#1}}
\makeatletter
\def\section{\def\@secnumfont{\mdseries}\@startsection{section}{1}%
  \z@{.7\linespacing\@plus\linespacing}{.5\linespacing}%
  {\normalfont\scshape\centering}}
\def\subsection{\def\@secnumfont{\bfseries}\@startsection{subsection}{2}%
  {\parindent}{.5\linespacing\@plus.7\linespacing}{-.5em}%
  {\normalfont\bfseries}}
\makeatother
\def\subl#1{\subsection{}\label{#1}}
 \nc{\Hom}{\operatorname{Hom}}
  \nc{\mode}{\operatorname{mod}}
\nc{\End}{\operatorname{End}} \nc{\wh}[1]{\widehat{#1}} \nc{\Ext}{\operatorname{Ext}} \nc{\ch}{\text{ch}} \nc{\ev}{\operatorname{ev}}
\nc{\Ob}{\operatorname{Ob}} \nc{\soc}{\operatorname{soc}} \nc{\rad}{\operatorname{rad}} \nc{\head}{\operatorname{head}}
\def\Im{\operatorname{Im}}
\def\gr{\operatorname{gr}}
\def\mult{\operatorname{mult}}
\def\Max{\operatorname{Max}}
\def\ann{\operatorname{Ann}}
\def\sym{\operatorname{sym}}
\def\loc{\operatorname{loc}}
\def\Res{\operatorname{\br^\lambda_A}}
\def\und{\underline}
\def\Lietg{$A_k(\lie{g})(\bsigma,r)$}
\def\res{\operatorname{res}}

 \nc{\Cal}{\cal} \nc{\Xp}[1]{X^+(#1)} \nc{\Xm}[1]{X^-(#1)}
\nc{\on}{\operatorname} \nc{\Z}{{\bold Z}} \nc{\J}{{\cal J}} \nc{\C}{{\bold C}} \nc{\Q}{{\bold Q}}
\renewcommand{\P}{{\cal P}}
\nc{\N}{{\Bbb N}} \nc\boa{\bold a} \nc\bob{\bold b} \nc\boc{\bold c} \nc\bod{\bold d} \nc\boe{\bold e} \nc\bof{\bold f} \nc\bog{\bold g}
\nc\boh{\bold h} \nc\boi{\bold i} \nc\boj{\bold j} \nc\bok{\bold k} \nc\bol{\bold l} \nc\bom{\bold m} \nc\bon{\bold n} \nc\boo{\bold o}
\nc\bop{\bold p} \nc\boq{\bold q} \nc\bor{\bold r} \nc\bos{\bold s} \nc\boT{\bold t} \nc\boF{\bold F} \nc\bou{\bold u} \nc\bov{\bold v}
\nc\bow{\bold w} \nc\boz{\bold z} \nc\boy{\bold y} \nc\ba{\bold A} \nc\bb{\bold B} \nc\bc{\mathbb C} \nc\bd{\bold D} \nc\be{\bold E} \nc\bg{\bold
G} \nc\bh{\bold H} \nc\bi{\bold I} \nc\bj{\bold J} \nc\bk{\bold K} \nc\bl{\bold L} \nc\bm{\bold M}  \nc\bo{\bold O} \nc\bp{\bold
P} \nc\bq{\bold Q} \nc\br{\bold R} \nc\bs{\bold S} \nc\bt{\bold T} \nc\bu{\bold U} \nc\bv{\bold V} \nc\bw{\bold W} \nc\bx{\bold
x} \nc\KR{\bold{KR}} \nc\rk{\bold{rk}} \nc\het{\text{ht }}
\nc\bz{\mathbb Z}
\nc\bn{\mathbb N}

\nc\toa{\tilde a} \nc\tob{\tilde b} \nc\toc{\tilde c} \nc\tod{\tilde d} \nc\toe{\tilde e} \nc\tof{\tilde f} \nc\tog{\tilde g} \nc\toh{\tilde h}
\nc\toi{\tilde i} \nc\toj{\tilde j} \nc\tok{\tilde k} \nc\tol{\tilde l} \nc\tom{\tilde m} \nc\ton{\tilde n} \nc\too{\tilde o} \nc\toq{\tilde q}
\nc\tor{\tilde r} \nc\tos{\tilde s} \nc\toT{\tilde t} \nc\tou{\tilde u} \nc\tov{\tilde v} \nc\tow{\tilde w} \nc\toz{\tilde z} \nc\woi{w_{\omega_i}}
\nc\chara{\operatorname{Char}}

\title[Demazure modules and Prime Representations]{Demazure modules of level two  and prime representations of  quantum affine $\lie{sl}_{n+1}$.}
\author[M. Brito, V. Chari, A. Moura]{Matheus Brito, Vyjayanthi Chari, and Adriano Moura}
\thanks{}

\address{Departamento de Matem\'atica, Unicamp, Campinas - SP - Brazil, 13083-859.}
\email{mbrito@ime.unicamp.br, aamoura@ime.unicamp.br}
\address{Department of Mathematics, University of California, Riverside, CA 92521, U.S.A.}
\email{vyjayanthi.chari@ucr.edu}
\thanks{M.B. was supported by FAPESP grant 2010/19458-9}
\thanks{V.C. was partially supported by DMS-1303052.}
\thanks{A.M. was partially supported by CNPq grant 303667/2011-7 and FAPESP grant 2014/09310-5}
\maketitle

\begin{abstract}
We study the  classical  limit of a  family of irreducible   representations of the quantum affine algebra associated to $\lie{sl}_{n+1}$.  After a suitable twist,  the  limit is    a module for $\lie{sl}_{n+1}[t]$, i.e., for  the maximal standard parabolic subalgebra of the  affine Lie algebra.  Our first result is about  the family of prime representations    introduced in \cite{hl:cluster},\cite{hl:cluster2}, in  the context of a monoidal categorification of cluster algebras. We show that these representations specialize (after twisting),  to  $\lie{sl}_{n+1}[t]$--stable, prime    Demazure modules in  level two integrable highest weight representations of the classical affine Lie algebra. It was proved in \cite{CSVW} that a  stable  Demazure module is isomorphic to the fusion product of stable, prime Demazure modules.  Our next result proves that such a fusion product   is the limit of the tensor product of the corresponding irreducible prime representations of quantum affine $\lie{sl}_{n+1}$.
\end{abstract}

\section*{Introduction}

%%%%%%%%%%%%INTRODUCTION%%%%%%%%%%%%%%%%%

The classification of  finite-dimensional irreducible representations of quantum affine algebras was  given in \cite{CPqa,CPbanff}. Since that time,  many different and  deep  approaches  have been developed to study these modules.  However, outside the simplest case of the quantum affine algebra associated to $\lie{sl}_2$, even the answer to a   basic question  such as a dimension formula,  is  not known  for an arbitrary irreducible representation. Thus, the focus of the study has been on particular families of modules: amongst the best known are the standard or local  Weyl modules,  Kirillov--Reshetikhin modules and their generalizations: the minimal affinizations (see \cite{Herkir,Heminaff,mukhinyoung,Nakajimakr} for instance). More recently, motivated by  categorification of cluster algebras, D. Hernandez and B. Leclerc identified   in \cite{hl:cluster} and \cite{hl:cluster2} an interesting class of \lq prime\rq\ irreducible  modules for the quantum affine algebra. (Recall that a representation is said to be prime if it is not isomorphic to a tensor product of non--trivial representations).  This class includes some of the known and well--studied prime representations but also included new families of examples.
One of the goals of  our  paper is to explore the structure of these prime representations by studying their $q\to 1$  specialization and to show that the character  of these modules is  given by the Demazure character formula.

The systematic study of the  classical limit of finite--dimensional representations of a quantum affine algebra was begun in  \cite{CPweyl}, where  a necessary and sufficient condition for the  existence of the limit was proved.  All  the interesting representations mentioned earlier  admit classical limits. The limit, when it exists,  is  a  finite--dimensional module for the corresponding affine Lie algebra.  Hence it is  also a module for the current algebra  $\lie g[t]$  of polynomial maps $\bc\to\lie g$, which is naturally a subalgebra of the affine  algebra.  (Here $\lie g$ is the underlying simple Lie algebra of the  affine Lie algebra). Equivalently  the current algebra is the commutator subalgebra of the standard maximal parabolic subalgebra of the affine Lie algebra.
The scaling element $d$ of the affine Lie algebra defines an integer grading on the affine Lie algebra and the current algebra is a graded subalgebra. The notion of the graded limit  of a representation of the quantum affine algebra was
developed in \cite{Cha01,CMkir}. It was shown in these papers that when $\lie g$ is of classical type, the Kirillov--Rehsetikhin modules could be regarded (after pulling  back by a suitable automorphism) as graded representations of the current algebra. The result is  now known in all types and we refer the reader to \cite{Kedem} for a detailed discussion.  	In \cite{Moura},  it  was shown in  some cases and conjectured in general  that minimal affinizations could also be regarded as a graded representation of the current algebra. The conjecture was established in \cite{Naoi2, Naoi3} when $\lie g$ is of classical type.

There is a well--known family of graded modules for the current algebra which arises as follows. Consider a highest weight irreducible integrable representation $V(\Lambda)$ of level $m$ for the affine Lie algebra.  Fix a Borel subalgebra of the affine Lie algebra. A  Demazure module of level $m$  is the module for the Borel subalgebra   generated by an extremal element of $V(\Lambda)$. Under natural conditions on the extremal vector, the Demazure module admits an action of the standard maximal parabolic subalgebra containing the chosen Borel subalgebra. Since $d$ is an element of the Borel subalgebra, it follows that such Demazure modules, which we call stable, are graded modules for $\lie g[t]$.

The relation between graded limits and Demazure modules was made in \cite{CL}. In that paper, it was shown, using results in \cite{CPweyl}, that any stable level one Demazure module is the graded limit of an irreducible local Weyl module (or  an irreducible standard module) for the quantum affine algebra. The results of \cite{FoL} imply that this remains true for simply--laced Lie algebras. In the non--simply laced case, this is no longer true; however it was shown in \cite{Naoi} that the graded limit admits a flag by level one Demazure modules.
The work of \cite{Cha01} and \cite{CMkir} shows that the graded limit of Kirillov--Reshetikhin modules, in the case when $\lie g$ is classical, is a higher level stable Demazure module. Finally, the work of \cite{Naoi2, Naoi3} shows that Demazure modules and generalized Demazure modules appear as graded limits of minimal affinizations.

We turn now to the new family of prime representations identified by Hernandez--Leclerc. The highest weight of such a representations satisfies a condition,  which is best described as minimal affinization by parts (see the next section for details).  In this paper we establish, in the case of $\lie{sl}_{n+1}$,  that the classical limit of such a prime representation is a stable level two prime Demazure module. The first hint of this connection comes from the work of \cite{CV} which gives a refined presentation of the stable Demazure modules.

We also answer in the affirmative the following question: does any level two Demazure module associated with affine $\lie{sl}_{n+1}$ appear as a classical limit of a (not necessarily prime)  representation of the quantum affine algebra?   This is a delicate question since the classical limit of a tensor product is not the tensor product of the classical limits.  The notion of a fusion product of $\lie g[t]$--modules was introduced in \cite{FL} and there are many examples which suggest that it  is  closely related to this  question.  In this paper we prove  that  certain  fusion products of  prime level two Demazure modules are  the classical limit of the tensor product of the corresponding representations of the quantum affine algebra.

A  detailed  overview of the results of this paper, a discussion of the  natural questions arising from our work and a description of the  overall organization of the paper is given in Section 1.

{\em Acknowledgement: This work was begun when the second and third authors were visiting   the Centre de Recherche Mathematique (CRM)  in connection with  the thematic semester \lq\lq New Directions in Lie theory\rq\rq. The authors are grateful to the CRM for the superb working environment.  }

\section{ The Main results} \label{mainresults} We describe  the main results of\ our paper and discuss  the  connections  with  \cite{CSVW}, \cite{hl:cluster} and \cite{hl:cluster2}.
Throughout the paper,  we denote by $\mathbb C$  the field of complex numbers, by $\mathbb Z$ the set of integers and by $\mathbb Z_+$,  $\mathbb N$ the set of non--negative and positive integers respectively.

\subsection{Simple, affine and current algebras}  \label{sln}
We shall only be interested in the  Lie algebra $\lie sl_{n+1}$ (denoted from now on as $\lie g)$   of $(n+1)\times (n+1)$ complex matrices of trace zero.   Let $I=\{1,\cdots, n\}$ be the set of vertices of the Dynkin diagram of $\lie g$ and $\{\alpha_i: i\in I\}$ and $\{\omega_i: i\in I\}$ be a set  of simple roots and fundamental weights.  The $\bz$ (resp. $\bz_+$) span of the simple roots will be denoted by $Q$ (resp. $Q^+)$  and the $\bz$ (resp. $\bz_+$) span of the fundamental weights is denoted by $P$ (resp. $P^+$). Define a partial order on $P$ by $\lambda\le \mu$ iff $\mu-\lambda\in Q^+$.  The positive roots are $$R^+=\{\alpha_{i,j}= \alpha_i+\alpha_{i+1}+ \cdots +\alpha_j: 1\le i\le j\le n\}.$$   Given   $\alpha_{i,j}=\alpha_i+\cdots+\alpha_j\in R^+$, let  $x^\pm_{ij}$ be the corresponding root vector of $\lie g$ and   set $x_i^\pm=x_{i,i}^\pm$ and   $h_i=[x_i^+, x_i^-]$. The elements $x_i^\pm, h_i$, $1\le i\le n$, generate $\lie g$ as Lie algebra.

 Let $\widehat {\lie g}$ be the untwisted affine Lie algebra associated to $\lie g$ which can be realized as follows.  Let $t$ be an indeterminate and let $\bc[t,t^{-1}]$ be the corresponding algebra of Laurent polynomials. Define a Lie algebra structure on the vector space $\lie g\otimes \bc[t,t^{-1}]\oplus\bc c\oplus \bc d $ by requiring $c$ to be central and setting $$[x\otimes t^r, y\otimes t^s]=[x,y]\otimes t^{r+s} +{\rm{tr}}(xy)c,\ \ [d,x\otimes t^r]=r(x\otimes t^r), \  \ x,y\in\lie g, \ \ r,s\in\bz.$$
The commutator subalgebra is $\lie g\otimes \bc[t,t^{-1}]\oplus\bc c $ and we shall denote it by $\widetilde{\lie g}$. We shall frequently regard the action of $d$  as defining a $\bz$--grading on $\widetilde{\lie g}$.

The current algebra  is the  $\bz_+$--graded  subalgebra $\lie g\otimes\bc[t]$ of $\widetilde{\lie g}$  and will be denoted as $\lie g[t]$.
We shall be interested in graded representations of  the current algebra: namely, $\bz$--graded vector spaces $V=\oplus_{s\in\bz} V[s]$ which admit an action of $\lie g[t]$ which is compatible with the grading,  $(\lie g\otimes t^r)V[s]\subset V[r+s]$, for all $r,s\in\bz$.  A map of graded  modules is a grade preserving map of $\lie g[t]$-modules.

\subsection{Quantized enveloping algebras and their $\ba$--forms}  Let  $\bc(q)$ be   the field of rational functions in an indeterminate $q$ and set  $\ba=\bz[q,q^{-1}]$. Let   $\bu_q(\lie g)$ and $\bu_q(\widetilde{\lie g})$ be the  quantized enveloping algebras  (defined over $\bc(q)$)  associated to $\lie g$ and $\widetilde{\lie g}$, respectively. The algebra $\bu_q(\lie g)$ is isomorphic to a subalgebra of $\bu_q(\widetilde{\lie g})$.     Let   $\bu_\ba(\lie g)$ and $\bu_\ba(\widetilde{\lie g})$ be the  $\ba$--form of  $\bu_q(\lie g)$ and $\bu_q(\widetilde{\lie g})$ defined in \cite{Lusztig}; these are free $\ba$--submodules such that $$\bu_q(\lie g)\cong \bu_\ba(\lie g)\otimes_{\ba}\bc(q),\qquad  \bu_q(\widetilde{\lie g})\cong \bu_\ba(\widetilde{\lie g})\otimes_{\ba}\bc(q).$$  Regard $\bc$ as an $\ba$--module by letting $q$ act as 1.   Then,
 $\bu_\ba(\widetilde{\lie g})\otimes_\ba\bc$ and $ \bu_\ba({\lie g})\otimes_\ba\bc$ are algebras over $\bc$ which have the universal enveloping algebra $\bu(\widetilde{\lie g})$ and  $\bu(\lie g)$ as canonical quotients.
Finally, recall that  $\bu_q(\widetilde{\lie g})$ is a Hopf algebra and that  $\bu_q(\lie g)$, $\bu_\ba(\widetilde{\lie g})$ and $\bu_\ba(\lie g)$ are all Hopf subalgebras.

\vskip 12pt

{\em Throughout the paper,  we shall as is  usual,  only be working with type one representations of quantized enveloping algberas and we will make no further mention of this fact}.

\subsection{ Finite--dimensional representations of $\lie g$, $\bu_q(\lie g)$ and their $\ba$--forms.}\hfill

It is well--known that the isomorphism classes of  irreducible finite--dimensional  representations of $\lie g$ and $\bu_q(\lie g)$ are indexed  by elements of $P^+$: given $\lambda\in P^+$, we denote by $V(\lambda)$ and $V_q(\lambda)$ an element of the corresponding isomorphism class.

It is also known that the category of finite--dimensional representations of these algebras is semisimple. Further,  the $\bu_q(\lie g)$--module $V_q(\lambda)$ admits an $\ba$--form $V_\ba(\lambda)$ which is a representation of $\bu_\ba(\lie g)$.   The space $V_\ba(\lambda)\otimes_{\ba}\bc$ is an irreducible  module for $\bu(\lie g)$ and we have, $$V_\ba(\lambda)\otimes_\ba\bc\cong_{\lie g} V(\lambda),\ \ \lambda\in P^+.$$

\subsection{The sets $\cal P^+\!\!\!$, $\cal P^+_\bz$, $\cal P^+_\bz(1)$  and the weight function.} \hfill

Let  $\cal P^+$ be the monoid consisting of  $n$--tuples of polynomials with coefficients in $\bc(q)[u]$ with constant term one and  with  coordinate-wise multiplication.   For $1\le i\le n$ and $a\in\bc(q)$, we take $\bomega_{i,a}$ to  be the $n$--tuple of polynomials where the only non--constant entry is the element $(1-au)$ in the $i$--th coordinate.  Let $\cal P^+_\bz$ be the submonoid of $\cal P^+$  generated by the elements $\bomega_{i,a}$, $1\le i\le n$ and $a\in q^{\bz}$.  Define $\wt:\cal P^+\to P^+$   by $\wt\bpi=\sum_{i=1}^n(\deg\pi_i) \omega_i$.
\begin{defn} Let $\cal P^+_\bz(1)$ be the subset of $\cal P^+_\bz$ consisting of the constant $n$--tuple and elements of the form  $\bomega_{i_1,a_1}\cdots\bomega_{i_k,a_k}$, where $1\le i_1<i_2<\cdots< i_k\le n$,  $a_j\in\bc(q)$,  $1\le k\le n$,  such that \begin{equation}\label{defcalp1}a_{i_j}a_{i_{j+1}}^{-1} =q^{\pm(i_{j+1}-i_j+2)},\ \ k\ge 2,\end{equation} and, if $j\le k-2$,  \begin{equation}\label{defcalp2}a_{i_j}a_{i_{j+1}}^{-1}= q^{\pm(i_{j+1}-i_j+2)}\implies a_{i_{j+1}}a_{i_{j+2}}^{-1}= q^{\mp(i_{j+2}-i_{j+1}+2)}.\end{equation}  \hfill\qedsymbol
\end{defn}
Note that \begin{equation}\label{p1}\wt\cal P^+_\bz(1)= P^+(1)=\{\lambda\in P^+: \lambda(h_i)\le 1, \ \ 1\le i\le n\}.\end{equation}

\subsection{Prime representations and prime factors.}  \hfill

 It was shown in \cite{CPqa}, \cite{CPbanff}, \cite{CPnato} that the isomorphism classes of irreducible finite--dimensional representations of $\bu_q(\widetilde{\lie g})$  are indexed by $\cal P^+$;
for  $\bpi\in\cal P^+$  let  $V(\bpi)$ be an element of the corresponding  isomorphism class. Note that the trivial representation corresponds to the constant $n$--tuple.
Given $\bpi,\bpi'\in\cal P^+$, the   tensor product $V(\bpi)\otimes V(\bpi')$ is generically irreducible and isomorphic to $V(\bpi\bpi')$. However,  necessary and sufficient conditions for this to hold are not known outside the case $n=1$ and motivated the interest in understanding the prime irreducible representations.
\begin{defn} We say that $V(\bpi)$ is a prime irreducible representation if it cannot be written in a non--trivial way as a tensor product of irreducible $\bu_q(\widetilde{\lie g})$ representations. We shall say that $\bpi_1,\cdots,\bpi_s$ are a set of prime factors of $\bpi$, if $V(\bpi_j)$ is prime and non--trivial, for all $1\le j\le s$, and  we have an isomorphism of $\bu_q(\widetilde{\lie g})$--modules
$$V(\bpi)\cong V(\bpi_1)\otimes\cdots\otimes V(\bpi_s).$$\hfill\qedsymbol
\end{defn}
Since $V(\bpi)$ is finite--dimensional it is clear that it is either prime or can be written as a tensor product of two or more non--trivial prime representations.  It is however, not clear that the set of prime factors of $V(\bpi)$ is unique: even in the case of simple Lie algebras, a unique factorization theorem for the tensor product of finite--dimensional irreducible modules was proved  relatively recently in \cite{Rajan}.

\subsection{Some  examples of prime representations.}\hfill

Regarding $V(\bpi)$ as a  finite--dimensional $\bu_q(\lie g)$--module, we can write  \begin{equation}\label{decomp} V(\bpi)\cong_{\bu_q(\lie{g})}V_q(\wt\bpi) \bigoplus_{\mu<\wt\bpi}\left(\dim\Hom_{\bu_q(\lie g)}(V_q(\mu),V(\bpi)\right) V_q(\mu).\end{equation} The best known examples of prime representations are the evaluation representations; namely an irreducible representation  $V(\bpi)$ of $\bu_q(\widetilde{\lie g})$ which is also irreducible as  an $\bu_q(\lie g)$--module, i.e., $V(\bpi)\cong_{\bu_q({\lie g})} V_q(\wt\bpi)$.   {\em It is important to recall here that we are in the case when $\lie g$ isomorphic to $\lie{sl}_{n+1}$ since our next assertion is false in the other types}. For  every $\lambda\in P^+$, there exist  elements  $\bpi\in\cal P^+$, $a\in\bc(q)$,  with $\wt\bpi=\lambda$  such that  $V(\bpi)\cong V_q(\lambda)$ as $\bu_q(\lie{g})$--modules.
An explicit formula for such elements was given in \cite{CPsmall} in complete generality. For the purposes of this paper, we will only need the following special cases.
\begin{lem}\label{minafformula}  Let $\bpi\in\cal P^+$.
\begin{enumerit}
\item[(i)] If  $\wt\bpi=\omega_i$, for some $1\le i\le n$, then $\bpi=\bomega_{i,a}$, for some $a\in\bc(q)$, and we have an isomorphism of $\bu_q(\lie{sl}_{n+1})$--modules,  $$V(\bpi)\cong V(\omega_i).$$
\item[(ii)] If $\wt\bpi=\omega_i+\omega_j$, for $1\le i\le j\le n$, then $$V(\bpi)\cong V_q(\omega_i+\omega_j)\iff
 \bpi= \bomega_{i,a}\bomega_{j, aq^{\pm(j-i+2)}},\ \ a\in\bc(q).$$\end{enumerit}In both cases, the module $V(\bpi)$ is prime.\hfill\qedsymbol
\end{lem}
\begin{rem} This Lemma is what motivates the definition of the set $\cal P_\bz^+(1)$.   Moreover, we shall see later that he module corresponding to an element of $V(\bpi)$, $\bpi\in\cal P_\bz^+(1)$ is \lq\lq minimal\rq\rq \ when restricted to suitable subalgebras of $\lie g$.
\end{rem}

\subsection{Further examples of  prime representations}\hfill

  Using Lemma \ref{minafformula}, it is easy to generate further examples of prime representations as follows.
  The next result was proved in   \cite{hl:cluster}, another proof can be found in \cite{cmy}.
\begin{lem} The module $V(\bpi)$ is prime for all $\bpi\in\cal P_\bz^+(1)$.\hfill\qedsymbol\end{lem}

Our main goal in this paper is to  understand  the $\bu_q(\lie g)$--character of the modules $V(\bpi)$, $\bpi\in\cal P_\bz^+(1)$, and, more generally, to relate it to other well--known modules for affine Lie algebras.

\subsection{The $\ba$--form $V_\ba(\bpi)$ and the representation $L(\bpi)$ of $\lie g[t]$.}\hfill

  It is in general not true  that an  arbitrary
irreducible finite--dimensional module for $\bu_q(\widetilde{\lie g})$ admits an $\ba$--form. In the special case, when $\bpi\in\cal P_\bz^+$, the results of
 \cite{Cha01}, \cite{CPweyl}, show  that $V(\bpi)$ does admit an $\ba$--form  and that $V_\ba(\bpi)\otimes_{\ba} \bc$ is an indecomposable and usually reducible module for  the enveloping algebra $\bu(\tilde{\lie g})$ and hence also for the current algebra $\lie g[t]$.  Consider  the  automorphism of  $\lie g[t]\to\lie g[t]$ defined  by  mapping $a\otimes t^r \to a\otimes (t-1)^r$, $a\in\lie g$, $r\in\bz_+$. Let $L(\bpi)$ be the representation of $\lie g[t]$ obtained by pulling back the $\lie g[t]$--module, $V_\ba(\bpi)\otimes_{\ba} \bc$ via this automorphism. Then, one can prove (see Section \ref{primereps}) that   $$(\lie g\otimes t^N\bc[t])L(\bpi)=0,\ \ \  N>>0.$$
Let $\delta_{r,s}$ be the usual Kronecker delta symbol. The following is  not hard to prove (see \cite{Moura} for instance).
\begin{lem} \label{grademinaff} Let $\bpi\in\cal P^+$ be such that $V(\bpi)\cong V_q(\wt\bpi)$. Then $$L(\bpi)\cong_{\lie g} V(\wt\bpi)\ \ {\rm and} \ \ (\lie g\otimes t\bc[t]) L(\bpi)=0.$$  In particular, $L(\bpi)$ is the graded $\lie g[t]$--module generated by an element $v_{\wt\bpi}$ of grade zero  and defining relations: \begin{gather*} (x_{1,n+1}^-\otimes t) v_{\wt\bpi}=0,\\ (x_i^+\otimes\bc[t])v_{\wt\bpi}=0,\ \ (h_i\otimes t^s)v_{\wt\bpi}=\delta_{s,0}(\deg\pi_i) v_{\wt\bpi},\ \ (x_i^-)^{\deg\pi_i+1} v_{\wt\bpi}=0,\end{gather*} for all $1\le i\le n$.\hfill\qedsymbol
\end{lem}

\subsection{The main result: a  presentation of $L(\bpi)$, $\bpi\in\cal P^+_\bz(1)$.}\label{dnulambda}\hfill

  Given $\nu\in P^+$ and $\lambda=\omega_{i_1}+\cdots +\omega_{i_k} \in\cal P^+(1)$, $1\leq i_1 < \cdots < i_k \leq n$,  let $M(\nu,\lambda)$ be the graded  $\lie g[t]$--module generated by an element $v_{\nu,\lambda}$ of grade zero with defining relations: for all $1\le i\le n$,  and $s\in\bz_+$, \begin{gather}\label{def1} (x^+_i\otimes\bc[t])v_{\nu,\lambda}=0,\ \ (h_i\otimes t^s)v_{\nu,\lambda}=\delta_{s,0}(\deg\pi_i) v_{\nu,\lambda},\ \  (x_i^-)^{\deg\pi_i+1} v_{\nu,\lambda}=0,\\ \label{def2}(x_i^-\otimes t)^{\nu(h_i)+\lambda(h_i)}v_{\nu,\mu}=0, \ \ 1\le i\le n \\ \label{def3}(x_{i_j,i_{j+1}}^-\otimes t^{\nu(h_{i_j}+h_{i_j+1}+\cdots+h_{i_{j+1}})+1})v_{\nu,\lambda}=0,\ \ 1\le j\le k-1.\end{gather}
In this paper, we shall prove,
\begin{thm}\label{first} Let $\bpi\in\cal P^+$ be such that $V(\bpi)$  has at most one prime factor in $\cal P^+_{\bz}(1)$ and all its other prime factors are of the form $\bomega_{i,a}\bomega_{i,aq^2}$, $1\le i\le n$, $a\in q^{\bz}$. Then $\wt\bpi=2\nu+\lambda$, for  a unique choice of $\nu\in P^+$ and $\lambda\in P^+(1)$, and we have an isomorphism of $\lie g[t]$--modules $$L(\bpi)\cong M(\nu,\lambda).$$  In particular, $L(\bpi)$ acquires the structure of a graded $\lie g[t]$--module.
\end{thm}

 \subsection{Demazure modules as fusion products.}\label{demazuremodules} \hfill

We assume for the moment that the reader is familiar with the notion of  $\lie g$--stable, Demazure modules  $D(\ell,\lambda)$  of level $\ell$, which  are indexed by pairs $(\ell,\lambda)\in\bn\times  P^+$. They  are  graded modules for $\lie g[t]$ and are generated by a vector $v_{\ell, \lambda}$:  a detailed development can be found in Section \ref{oddtpeven} where we shall also see that  we have an isomorphism of $\lie{sl}_{n+1}$--modules,  $$D(\ell, \ell\omega_i)\cong V(\ell\omega_i), \ \  1\le i\le n,\ \  \ell\in\bn.$$      The tensor product of level $\ell$--Demazure modules is not a level $\ell$--module. In \cite{FL}, the authors introduced a new $\lie g[t]$--structure on the tensor product of cyclic graded $\lie g[t]$--modules; the resulting $\lie g[t]$--module (unlike the tensor product) is a {\em cyclic} $\lie g[t]$--module. This structure, called the fusion product,  depends on a choice of  complex numbers, a distinct one for each factor in the tensor product.   The underlying $\lie g$--module structure is  unchanged and in  many cases, it is known that the fusion product is independent of this choice of parameters; the case of interest to us is the following special case of a  result  proved in  \cite{CSVW}.
\begin{prop}\label{level2fusion} Let $\nu=\sum_{i=1}^nr_i\omega_i\in P^+$ and $\lambda\in P^+(1)$. Then $$D(2,2\nu+\lambda)\cong D(2,2\omega_1)^{*r_1}*\cdots *D(2,2\omega_n)^{*r_n}*D(2,\lambda).$$Moreover, the module $D(2,\lambda)$ is prime in the sense that it is not isomorphic to a tensor product of non--trivial $\lie g$--modules. \hfill\qedsymbol
\end{prop}
\subsection{The connection with Demazure modules.}\hfill

 The connection between Theorem \ref{first} and   Demazure modules is made via the following proposition  proved in Section \ref{effpres} of this paper.
\begin{prop}\label{second} For  $\nu\in P^+$ and $\lambda\in P^+(1)$, we have an isomorphism $$M(\nu,\lambda)\cong_{\lie g[t]} D(2, 2\nu+\lambda).$$\end{prop}

\subsection{The connection with the category $\cal C_\kappa$. }\hfill

 We discuss the relationship of our work with that of  \cite{hl:cluster, hl:cluster2}.
Let  $\kappa: \{1,2,\cdots, n\}\to \bz$ be a function  satisfying $|\kappa(i+1)-\kappa(i)|= 1$, for $1\le i\le n-1$.   Let $\cal C_\kappa$ be  the full subcategory of  finite--dimensional representations  of $\bu_q(\widetilde{\lie g})$ defined as follows: an object of $\cal C_\kappa$   has all its Jordan--Holder components of the form $V(\bpi)$, where  $\bpi\in\cal P^{+}$ is a  product of terms of   the form  $\bomega_{i,a}$,  $a\in\{q^{\kappa(i)}, q^{\kappa(i)+2}\}$, $1\le i\le n$..  The following is  a slight reformulation of  results proved in  \cite{hl:cluster, hl:hall, hl:cluster2}.\begin{thm}\label{hl}  The category $\cal C_\kappa$ is closed under taking tensor products. An irreducible object of $\cal  C_\kappa$ is a tensor product of prime irreducible objects of $\cal C_\kappa$.  An irreducible object $V(\bpi)$ of  $\cal C_\kappa$  is prime only if  $\bpi\in\cal P_\bz^+(1)$ or if $\bpi=\bomega_{i,q^{\kappa(i)}}\bomega_{i, q^{\kappa(i)+2}}$, for $1\le i\le n$. Moreover, given $\bpi\in\cal P^+_\bz(1)$, there exists a height function $\kappa$ such that $V(\bpi)$ is a prime object of $\cal C_\kappa$. \hfill\qedsymbol
\end{thm}
We make some comments about our reformulation.  In  \cite{hl:cluster2}, the authors define a quiver $\mathbf Q_\kappa$ whose vertices  are the elements of the set $\{1,\cdots, n\}$ and the edges are the set $i\to i+1$, if $\kappa(i)<\kappa(i+1)$, and $i+1\to i,$ otherwise. Given any subset $J=(1\le i_1<i_2<\cdots <i_k\le n)$ of $I$, consider the connected  subquiver determined by this subset and let $J_<$ and $J_>$ be the set of sinks and sources respectively, of this subquiver. According to \cite{hl:cluster2}, the representation $V(\bpi)$ associated to $$\bpi=\prod_{j\in J_<}\bomega_{j,q^{\kappa (j)}}\prod_{j\in J_>}\bomega_{j,q^{\kappa(j)+2}},$$ is prime and, moreover, any prime object in $\cal C_\kappa$ is either of this form or is isomorphic to $V(\bomega_{i,q^{\kappa(i)}}\bomega_{i, q^{\kappa(i)+2}})$, for $1\le i\le n$. It is straightforward to check that the element $\bpi\in\cal P^+_\bz(1)$. Conversely, given $\bpi=\bomega_{i_1,a_1}\cdots\bomega_{i_k,a_k}\in\cal P_\bz^+(1)$, consider the height function $\kappa$  given by requiring the elements $i_1,i_3,\cdots$ to be the sinks of $\mathbf Q_\kappa$ and $i_2, i_4,\cdots$ to be the sources of the quiver. Then $V(\bpi)$ is a prime object of $\cal C_\kappa$.

Very little has been  known  so far about prime objects $V(\bpi)$ in $\cal C_\kappa$, except in the case when $\bpi=\bomega_{i,\kappa(i)}\bomega_{j, \kappa(j)}$ or $\bpi=\bomega_{i,q^{\kappa(i)}}\bomega_{i, q^{\kappa(i)+2}}$, for $1\le i\le n$, where we can use Lemma \ref{minafformula}.  As a consequence of    Lemma \ref{grademinaff} and  Theorem  \ref{first} our results  show that if $V(\bpi)$ is a prime object of $\cal C_\kappa$ then one knows a presentation of  $L(\bpi)$ and that it is a graded  $\lie g[t]$--module.  Moreover,  since
the $\bu_q(\lie g)$-character of $V(\bpi)$ is  the same as the $\lie g$--character of $L(\bpi)$, Propositon \ref{second} shows that  their character is given by  the Demazure character formula. Together with Proposition \ref{level2fusion}, we see that $V(\bpi)$, $\bpi\in\cal P^+_\bz(1)$, is  {\em strongly prime},  in the sense that $V(\bpi)$ is not isomorphic to a tensor product of non--trivial $\bu_q(\lie g)$--modules. To the best of our knowledge all known examples of prime representations  are strongly prime. But,  it is far from clear  that the two notions are equivalent.

\subsection{An outline of the proof.}\label{outline}\hfill

  We outline the main steps of the proof of Theorem \ref{first}. It generalizes   ideas in \cite{Moura} (see also  \cite{Naoi2}) where a similar question was studied for a different family of irreducible representations.

 Suppose that $\bpi_1,\bpi_2\in\cal P_\bz^+$ are such that we have an injective map of $\bu_q(\widetilde{\lie g})$-modules $V(\bpi_1\bpi_2)\to V(\bpi_1)\otimes V(\bpi_2)$. Since $\bu_\ba(\widetilde{\lie g})$ is a Hopf subalgebra of $\bu_q(\widetilde{\lie g})$, we get an injective map $V_\ba(\bpi_1\bpi_2)\to V_\ba(\bpi_1)\otimes V_\ba(\bpi_2)$.  It was shown in \cite[Lemma 2.20, Proposition 3.21]{Moura} that tensoring with $\otimes_\ba\bc$ and pulling back by the automorphism of $\lie g[t]$ induced by  $t\to t-1$,   gives rise to a map of $\lie g[t]$--modules $L(\bpi)\to L(\bpi_1)\otimes L(\bpi_2)$. This map is neither injective nor surjective even  when $V(\bpi)\cong V(\bpi_1)\otimes V(\bpi_2)$, but  plays   big role in this paper.  The main steps in the proof of Theorem \ref{first} are the following. Let  $\bpi\in\cal P_{\bz}^+(1)$.  \\

(i) There exists a surjective map of $\lie g[t]$--modules $\varphi_1: M(0,\wt\bpi)\to L(\bpi)\to 0$.\\\

(ii)   There exists  $\bpi^o,\bpi^e\in\cal P_\bz^+(1)$  with  $\bpi=\bpi^o\bpi^e$ such that we have
an injective map of $\bu_q(\widetilde{\lie g})$--modules $V(\bpi)\to V(\bpi^o)\otimes V(\bpi^e)$. Moreover,
$$L(\bpi^o)\cong_{\lie g[t]}  D(1,\wt\bpi^o),\qquad\ L(\bpi^e)\cong_{\lie g[t]} D(1,\wt\bpi^e).$$
The  induced map $$\varphi_2: L(\bpi)\to L(\bpi^0)\otimes L(\bpi^e)\cong D(1,\wt\bpi^o)\otimes D(1,\wt\bpi^e),$$ satisfies, $$\varphi_2(\varphi_1(v_{\wt\bpi}))= v_{1,\wt\bpi^o}\otimes v_{1,\wt\bpi^e}.$$

(iii) There exists an injective map of $\lie g[t]$--modules $$D(2,\wt\bpi)\to D(1,\wt\bpi^o)\otimes D(1,\wt\bpi^e),\ \ {\rm{with}}\ \ v_{2,\wt\bpi}\to v_{1,\wt\bpi^o}\otimes v_{1,\wt\bpi^e}.$$

Consider the composite map $$\varphi_2\varphi_1: M(0,\wt\bpi)\to  D(1,\wt\bpi^o)\otimes D(1,\wt\bpi^e).$$  The preceding steps show that the image of this map is $D(2,\wt\bpi)$.
 Proposition \ref{second} proves that $\varphi_2\varphi_1$ must be injective and hence it follows that $\varphi_1$ is  an isomorphism  proving Theorem \ref{first} when $\bpi\in\cal  P_\bz^+(1)$.\\

(iv) The next  step is to prove Proposition \ref{second}.\\

(v) The final step is  to  prove  the following.  Suppose that $V$ is a cyclic $\lie g[t]$-module generated by a vector $v$ satisfying \eqref{def1}. Assume that $V\cong D(2,2\nu+\lambda)$ as $\lie g$--modules. Then, $V\cong D(2,2\nu+\lambda)$ as $\lie g[t]$--modules.\\

We now  deduce Theorem \ref{first} in full generality. We prove in Lemma \ref{lpi} that $L(\bpi)$ is a cyclic $\lie g[t]$--module.
The assumptions on $\bpi$ in Theorem \ref{first} imply that $V(\bpi)$ is isomorphic, as  $\bu_q(\widetilde{\lie g})$--modules, to a tensor product of modules of the form $V(\bomega_{i,a}\bomega_{i,aq^2})$, $1\le i\le n$, along with a single module $V(\bpi_1)$, with $\bpi_1\in\cal P^+_\bz(1)$. Hence, we can write $\wt\bpi=2\nu+\wt\bpi_1$, with $\nu=\sum_{i=1}^nr_i\omega_i\in P^+$ and  $\wt\bpi_1\in P^+(1)$.
Since this is an isomorphism also of $\bu_q(\lie g)$--modules, it follows that we have an isomorphism of $\lie g$ (but not of $\lie g[t])$) modules, $$L(\bpi)\cong_{\lie g} L(2\omega_1)^{\otimes r_1}\otimes\cdots\otimes L(2\omega_n)^{\otimes r_n}\otimes L(\bpi_1).$$
 Using Lemma \ref{minafformula}, the fact that $D(2,2\omega_i)\cong_{\lie g}V(2\omega_i)$ and  Theorem \ref{first} for $\bpi_1$,  we get  $$L(\bpi) \cong_{\lie g} D(2, 2\omega_1)^{\otimes r_1}\otimes \cdots\otimes D(2,2\omega_n)^{\otimes r_n}\otimes D(2, \wt\bpi_1).$$ Using Proposition \ref{level2fusion} gives $$L(\bpi)\cong_{\lie g} D(2, \wt\bpi).$$
Using step (iv) now proves that $L(\bpi)$ is isomorphic to $D(2,\wt\bpi)$ as $\lie g[t]$--modules and the proof is complete.
\\

The proof of steps (i) and (ii) are in Section \ref{primereps}, the proof of step (iii) is in Section \ref{oddtpeven}, the proof of step (iv) is in Section  \ref{effpres} and the the proof of step (v) is in Section \ref{fusionproducts}.

 \subsection{Further Questions.} \hfill

As we have remarked, the isomorphism  $$V(\bpi_1\bpi_2)\cong V(\bpi_1)\otimes V(\bpi_2),$$ does not imply that $L(\bpi_1\bpi_2)$ is isomorphic $ L(\bpi_1)\otimes L(\bpi_2)$.  However,  in the case when $L(\bpi_j)$ are graded $\lie g[t]$--modules, it is true  in many known examples,  that $$L(\bpi_1\bpi_2)\cong L(\bpi_1)*L(\bpi_2).$$ Theorem \ref{first}, Proposition \ref{level2fusion}  and  Proposition  \ref{second} add to  these growing number of examples which support the conjecture that such a statement might be  true in general.

We now discuss some natural questions arising from our work.  The Demazure character formula is not easy to compute and an interesting question would be to determine the $\lie g$--module decomposition of the prime  level two Demazure modules. Preliminary calculations show that there could be some interesting combinatorics associated with it, along the lines of the formulae given for the well--known Kirillov--Reshetikhin modules.

Another direction is the following: there exists irreducible (but not prime) objects of $\cal C_\kappa$ which are not level two Demazure modules.   In fact any tensor product of modules which has more than one prime factor of the form $V(\bpi)$, $\bpi\in\cal P^+_\bz(1)$,  does not usually specialize to a  level two  Demazure module. It is natural to speculate that they too have graded limits which are also relates in some way to other known representations for the current algebra.  In simple cases, it appears that  these are    quotients of  a family  modules for $\lie{sl}_{n+1}[t]$ defined and studied in \cite{CV}; they also appear to be related to the generalized Demazure modules studied in \cite{LLM} and \cite{Naoi}.

There are two more very obvious questions: can one formulate and prove analogous results for Demazure modules of level at least three for $\lie{sl}_{n+1}[t]$ and are there analogous results for the other simple Lie algebras? It is most convenient to address these two questions together.  Assume therefore that $\lie g$ is an arbitrary simple Lie algebra.  If  $\lie g$  is simply--laced  then one uses \cite{CPweyl}, \cite{CL} and \cite{FoL}, to prove that ,   any level one Demazure module is isomorphic to the graded limit  of a module for the quantum affine algebra. However, in the non--simply laced case this is false; the level one Demazure modules are too small. It was shown  in \cite{Naoi} that the graded limit usually admits only a flag by Demazure modules.

This phenomenon persists as we move to higher levels. Thus for type $D$ and $E$, the level two Demazure modules are again too small; one sees this already in the case of most of  the Kirillov-Reshetikhin modules.  Similarly for type $A$ beyond level two, some of  the Demazure modules   appear to be  too small. In all these cases there is evidence to suggest that the  generalized Demazure modules and the modules defined in \cite{CV} might  be the correct objects. In particular, in work in progress (see also \cite{britotheses})  we show that the minimal affinizations, which are known through the work of \cite{Naoi} to specialize to generalized Demazure module are isomorphic to the modules defined in \cite{CV}.

\section{ Prime representations and graded limits}\label{primereps} In this section we recall the necessary results from the theory of finite--dimensional representation of quantum affine $\lie{sl}_{n+1}$.  The results can  be stated without introducing in too much detail,  the extensive notation of the quantum affine algebras. We  refer the reader to  \cite{CPbook} for the basic definitions and  to \cite[Section 3] {Naoi2} for an excellent exposition with detailed references, for the results discussed here,  on graded limits and  minimal affinizations.

\subsection{The modules $L(\bpi)$ for  $\bpi\in\cal P^+_\bz$.}\label{lpi} \hfill

We begin the section by elaborating on the definition of the $\lie g[t]$--modules $L(\bpi)$, $\bpi\in\cal P^+_\bz$. Special cases of what we are going to say are in the literature (\cite{Cha01, Moura, Naoi}) but  not  in the generality we need.

For $\bpi\in \cal P^+$, let $\bar\bpi=(\pi_1(u),\cdots \pi_n(u))\in\bc[u]$ be the $n$--tuple of polynomials with complex coefficients obtained from $\bpi$ by setting $q=1$. If  $\bpi\in\cal P_\bz^+$, then it follows that $\bar\pi_i(u)=(1-u)^{\deg\pi_i}$, for all $1\leq i \leq n$.
 It was shown in \cite[Section 4]{CPweyl} that the module $V(\bpi)$, $\bpi\in\cal P^+_\bz$, has an $\ba$--form $V_\ba(\bpi)$ which gives  an action of $\widetilde{\lie g}$ on $\overline {V(\bpi)}:=V_\ba(\bpi)\otimes_\ba \bc$.
Moreover, it was proved that $\overline {V(\bpi)}$ is generated, as a $\widetilde{\lie g}$--module, by an element $v_{\overline{\bpi}}$ which satisfies the relations:  $$(x_i^+ \otimes\bc[t,t^{-1}])v_{\overline{\bpi}}=0,\ \  (h_i\otimes t^r)v_{\overline{\bpi}} =  (\deg\pi_i) v_{\overline{\bpi}},\ \ (x_i^-\otimes 1)^{\deg\pi_i+1}v_{\overline{\bpi}}=0,$$  for all $1\le i\le n$. Since $\overline {V(\bpi)}$ is a finite--dimensional module for $\widetilde{\lie g}$ it follows  that the central element acts as zero, i.e., $\overline {V(\bpi)}$  is a module for $\lie g\otimes\bc[t,t^{-1}]$.
 Hence, there must exist $f\in\bc[t,t^{-1}]$ of minimal degree such that $$(\lie g\otimes f\bc[t,t^ {-1}]) \overline {V(\bpi)}=0.$$  Proposition 2.7 of \cite{CFS} shows that  $f$ must be of the form $(t-1)^N$ for some $N>>0$.
 This means  that the module $\overline{V(\bpi)}$ is a module for the quotient Lie algebra $\lie g\otimes \bc[t,t^{-1}]/(t-1)^N$. Using the isomorphism, $$\lie g\otimes \bc[t,t^{-1}]/(t-1)^N\cong\lie g\otimes \bc[t]/(t-1)^N,$$ we see that we can regard $\overline{V(\bpi)}$ as a module for $\lie g[t]$ generated by the element $v_{\overline\bpi}$
 such that $$(\lie g\otimes (t-1)^N\bc[t]) \overline {V(\bpi)}=0.$$  Pulling back  $\overline{V(\bpi)}$ via the automorphism $x\otimes t^r\to x\otimes (t-1)^r$ of $\lie g[t]$, we get the module which we have called  $L(\bpi)$. Summarizing we have,
\begin{lem} The $\lie g[t]$--module $L(\bpi)$ is generated by an element $v_{\overline{\bpi}}$ satisfying   $$(x_i^+ \otimes\bc[t])v_{\overline{\bpi}}=0,\ \ \ \   (h_i\otimes t^r)v_{\overline{\bpi}} =  \delta_{r,0}(\deg\pi_i)  v_{\overline{\bpi}},\ \ \ \ (x_i^-\otimes 1)^{\deg\pi_i+1}v_{\overline{\bpi}}=0,$$ for all $1\le i\le n$.  In particular,  we have $$ L(\bpi)_{\wt\bpi}=\{v\in L(\bpi): (h\otimes 1)v=\wt\bpi(h)v, \ \ h\in \lie h\}= \bc v_{\overline{\bpi}}.$$\hfill\qedsymbol
\end{lem}

\subsection{Restrictions to subalgebras} \label{diagramsubalgebras} \hfill

 Suppose that $J$ is a connected subset of $I $ of cardinality $m$. Then the subalgebra of $\lie g$ generated by the elements $x_j^\pm$, $j\in J$, is isomorphic to $\lie{sl}_{m+1}$ and  will be denoted by $\lie g_J$.  We have canonical inclusions of the quantized enveloping algebras, $\bu_q(\lie g_J)\hookrightarrow\bu_q(\lie g)$ and $\bu_q(\widetilde{\lie g_J})\hookrightarrow\bu_q(\widetilde{\lie g})$.

The $\bz_+$--span of the weights $\omega_j$, $j\in J$, is denoted by  $P_J^+$; the subsets $Q_J^+$ of $Q^+$ and $\cal P_J^+$ of $\cal P^+$ are defined in the obvious way.   Define surjective maps $P^+\to P_J^+$ and $\cal P^+\to\cal P^+_J$ by
$$\lambda=\sum_{i=1}^nr_i\omega_i\to\lambda_J=\sum_{j\in J}r_j\omega_j,\qquad \bpi=(\pi_1,\cdots,\pi_n)\to \bpi_J=(\pi_j)_{j\in J}.$$
 The set $P_J^+$ (resp. $\cal P_J^+$)   indexes the isomorphism classes of finite--dimensional  irreducible representations of $\lie g_J$ and $\bu_q(\lie g_J)$ (resp.  of $\bu_q(\widetilde{\lie g_J})$). Given $\lambda\in P^+_J$, let $V^J(\lambda)$ be a finite--dimensional irreducible representation of $\lie g_J$; the modules $V^J_q(\lambda)$ and $V^J(\bpi)$, $\bpi\in\cal P^+_J$, are defined in the obvious way.
The following result is proved by standard methods (see \cite{CPminaff} for details in the quantum case).
\begin{prop} Let $\lambda\in P^+$ and $\bpi\in\cal P^+$,
\begin{enumerit}
\item[(i)] We have an inclusion of $\lie g_J$--modules $V^J(\lambda_J)\to V(\lambda)$. Analogous statements hold for $V_q^J(\lambda_J)$ and for $V^J(\bpi_J)$.
\item[(ii)] Suppose that $\mu\in P^+$ is such that $\lambda-\mu\in Q^+_J$. Then $$\Hom_{\bu_q({\lie g})}(V_q(\mu), V(\bpi))\ne 0\implies \Hom_{\bu_q{(\lie g_J})}(V_q^J(\mu_J), V^J(\bpi_J))\ne 0.$$
\end{enumerit}\hfill\qedsymbol
\end{prop}

	\subsection{Existence of $\varphi_1$}\hfill

 Suppose now that $\bpi_1\in\cal P^+_\bz(1)$ and let $\lambda:=\wt\bpi=\omega_{i_1}+\cdots +\omega_{i_k}$, for some $1\le k\le n$ such that $1\leq i_1< \cdots < i_k\leq n$.  In view of the definition of $M(0,\wt\bpi)$ given in Section \ref{dnulambda} and  Lemma  \ref{lpi}, the existence of the surjective map  $\varphi_1: M(0,\wt\bpi)\to L(\bpi)$ follows, if  we prove that $$(x_i^-\otimes t^{\lambda(h_{i})}) v_{\overline{\bpi}}=0,\ \ \  1\le i\le n, \ \
(x^-_{i_j,i_{j+1}}\otimes t)v_{\overline{\bpi}}=0,\ \ 1\le j\le k-1.$$  If $\lambda(h_i)=0$ then the first equality is clear from Lemma \ref{lpi}. If $\lambda(h_i)=1$, i.e., $i=i_j$ for some $1\le j\le k$, then Lemma \ref{lpi} gives $$(x_{i}^-\otimes 1)^2v_{\overline{\bpi}}=0.$$ Applying $(x_i^+\otimes t)$, taking commutators and using the relations in Lemma \ref{lpi} again proves that $(x_i^-\otimes t) v_{\overline{\bpi}}=0$.

Suppose for a contradiction that $v:= (x^-_{i_j,i_{j+1}}\otimes t)v_{\overline{\bpi}}\ne 0$, for some $1\le j\le k-1$. Then, we claim that $x_s^+v=0,\ 1\le s\le n.$ This is obviously true except when $s=i_j$ or $s=i_{j+1}$, when we have to prove that$$(x^-_{i_j+1,i_{j+1}}\otimes t)v_{\overline{\bpi}}= 0= (x^-_{i_j,i_{{j+1}}-1}\otimes t)v_{\overline{\bpi}}.$$ If $i_{j+1}=i_j+1$, i.e., $x^-_{i_j+1,i_{j+1}}=x^-_{i_{j+1}}$ or $x^-_{i_j,i_{j+1}-1} =x^-_{i_j}$,  this was  established in the previous paragraph.

If $i_{j+1}>i_j+1$ we can write $$(x^-_{i_j+1,i_{j+1}}\otimes t)= [x^-_{i_j+1,i_{j+1}-1}, x^-_{i_{j+1}}\otimes t].$$ Since $\deg\pi_{i}=0$ if $i_j<i<i_{j+1}$, we have $x^-_{i_j+1,i_{j+1}-1}v_{\overline{\bpi}}=0$. Together with the fact that we have proved $(x^-_{i_{j+1}}\otimes t)v_{\overline{\bpi}}=0$, it follows that $(x^-_{i_j+1,i_{j+1}}\otimes t)v_{\overline{\bpi}}=0$.
%If $i_{j+1}>i_j+1$ we can write $$(x^-_{i_j+1,i_{j+1}}\otimes t)=  [x^-_{i_j+1}, x^-_{i_j+2,i_{j+1}}\otimes t].$$  Since $\deg\pi_{i_{j}+1}=0$ we have $x^-_{i_j+1}v_{\overline{\bpi}}=0$ and we have proved $(x^-_{,i_{j+1}}\otimes t)v_{\overline{\bpi}}=0$ it follows that $(x^-_{i_j+1,i_{j+1}}\otimes t)v_{\overline{\bpi}}=0$.
The proof that $ (x^-_{i_j,i_{{j+1}}-1}\otimes t)v_{\overline{\bpi}}=0$ is similar and we omit the details. This proves the claim and
it  follows that we have $$\dim\Hom_{\lie g}(V(\mu), L(\bpi)) >0,\ \ \  \mu=\lambda-\sum_{s=0}^{i_{j+1}-i_j}\alpha_{i_j+s},$$ and hence, also \begin{equation}\label{imp}\dim\Hom_{\bu_q(\lie g)}(V_q(\mu), V(\bpi)) >0,\ \ \mu=\lambda-\sum_{s=0}^{i_{j+1}-i_j}\alpha_{i_j+s}.\end{equation}
 Setting $J=\{i_j, i_{j}+1,\cdots, i_{j+1}\}$, we have $$\mu_J=0,\qquad \bpi_J=(\bomega_{i_j, a_j},1,\cdots, 1, \bomega_{i_{j+1}, a_{j+1}}), \ \ a_{i_j}a_{i_{j+1}}^{-1}= q^{\pm( i_{j+1}-i_j+2)}.$$  Using Proposition \ref{diagramsubalgebras} and equation  \eqref{imp} we have \begin{equation}\label{imp1}\dim\Hom_{\bu_q(\lie g_J)}(V^J_q(0), V^J(\bpi_J)) >0.\end{equation} On the other hand,  Lemma \ref{minafformula} applies to the representation  $V^J(\bpi_J)$ of $\bu_q(\widetilde{\lie g}_J)$,   and, hence,  we have $V^J(\bpi_J)\cong V^J_q(\lambda_J)$ as  $\bu_q(\lie g_J)$--modules. This contradicts equation \eqref{imp1} and hence  we must have $v=0$.  The existence of the surjective map $\varphi_1:M(0,\wt\bpi)\to L(\bpi)\to 0$ is established.

\subsection{Tensor products and defining relations}  \hfill

The relations given in Lemma \ref{lpi} are not necessarily the  defining relations of $L(\bpi)$. However, we have the following result, which is a  special case of \cite[Section 4]{CPweyl} and the main result of \cite{CL}.
\begin{thm}\label{primefund}  Let $\bpi\in\cal P_\bz^+$ be such that the prime factors  of $V(\bpi)$ are $$\{V(\bomega_{j_s,b_s}): 1\le s\le m, \ \ 1\le j_s\le n,\ \  b_s\in\ q^{\bz}\},$$ i.e.,  $$V(\bpi)\cong_{\bu_q(\widetilde{\lie g})} V(\bomega_{j_1,b_1})\otimes\cdots \otimes V(\omega_{j_m, b_m}).$$  Then $L(\bpi)$ is generated by the element $v_{\overline{\bpi}}$ with {\em{defining}} relations:$$(x_i^+ \otimes\bc[t])v_{\overline{\bpi}}=0,\ \ \ \   (h_i\otimes t^r)v_{\overline{\bpi}} =  \delta_{r,0}(\deg\pi_i)  v_{\overline{\bpi}},\ \ \ \ (x_i^-\otimes 1)^{\deg\pi_i+1}v_{\overline{\bpi}}=0,$$ for all $1\le i\le n$. \hfill\qedsymbol
 \end{thm}
\subsection{Simple socle in a tensor product}\hfill

 The hypothesis  in Theorem \ref{primefund}  that $V(\bpi)$ has all its prime factors of the form $\bomega_{i,a}$,  $1\le i\le n$, $a\in q^{\bz}$, is generically  true.  This is  seen from the next result,   the dual of which is proved in \cite[Theorem 3 and  Corollary 5.1]{CBraid}.
\begin{prop}\label{cbraidthm}  Let $m\ge 1$,  $1\le j_1,\cdots, j_m\le n$ and $b_1,\cdots, b_m \in q^{\bz}$, be such that\\  \begin{equation}\label{singular}  s>r\implies   b_r/b_s\notin{\{ q^{2p+2-j_s-j_r}: \max\{j_r,j_s\} <p+1\le \min\{j_r+j_s, n+1\}\}}.\end{equation}\\
Then $V(\bomega_{j_1,b_1}\cdots\bomega_{j_m,b_m})$ is the unique irreducible submodule of $ V(\bomega_{j_1,b_1})\otimes\cdots\otimes V(\bomega_{j_m,b_m})$.
  Moreover, if equation \eqref{singular} holds for all $1\le r,s\le m$, then we have an isomorphism of $\bu_q(\widetilde{\lie g})$--modules,$$V(\bomega_{j_1,b_1})\otimes\cdots\otimes V(\bomega_{j_m,b_m})\cong V(\bpi).$$

\hfill\qedsymbol
\end{prop}

\subsection{The elements $\bpi^o$ and $\bpi^e$}\hfill

  Let $\bpi\in\cal P^+_{\mathbb Z}(1)$, in which case we can write,   \begin{gather*}\bpi=\bomega_{i_1, a_1}\cdots \bomega_{i_k,a_k},\quad
 1\le i_1<\cdots <i_k\le n,\end{gather*} where either \begin{gather*}  a_{1} = q^m,\ \ \ \   {a_{2j}}/{a_{2j-1}}  =q^{i_{2j}-i_{2j-1}+2},
 \ \ \ \  {a_{2j+1}}/{a_{2j}} =q^{i_{2j}-i_{2j+1}-2},\end{gather*} or\begin{gather*}  a_{1} = q^m,\ \ \ \   {a_{2j}}/{a_{2j-1}}  =q^{-(i_{2j}-i_{2j-1}+2)},
 \ \ \ \  {a_{2j+1}}/{a_{2j}} =q^{-(i_{2j}-i_{2j+1}-2)},\end{gather*} for some $m\in \bz$. We shall always assume that we are in the first case. The proof in the other case is identical.

Define  $$\bpi^o=\begin{cases}\bomega_{i_1,a_1}\bomega_{i_3,a_3}\cdots \bomega_{i_k,a_k},\ \ k \ {\rm{odd}}\\
\bomega_{i_1,a_1}\bomega_{i_3,a_3}\cdots \bomega_{i_{k-1},a_{k-1}},\ \   k \ {\rm{even}}\end{cases},\qquad \bpi^e=\begin{cases}\bomega_{i_2,a_2}\bomega_{i_4,a_4}\cdots \bomega_{i_{k-1},a_{k-1}},\ \ k \ {\rm{odd}}\\
\bomega_{i_2,a_2}\bomega_{i_4,a_4}\cdots \bomega_{i_k,a_k},\ \   k \ {\rm{even}}\end{cases},$$\\ and note that $\bpi^e\bpi^o=\bpi$.

\begin{prop}\label{cyclic}  Assume that $k$ is even. We have an isomorphism of $\bu_q(\widetilde{\lie{g}})$--modules,  $$ V(\bpi^o) = V(\bomega_{i_1, a_1})\otimes \cdots\otimes V(\bomega_{i_{k-1},a_{k-1}}),\qquad V(\bpi^e)\cong V_q(\bomega_{i_2,a_2})\otimes\cdots\otimes V(\bomega_{i_k,a_k}).$$ Moreover $V(\bpi)$ is a submodule  of $V(\bpi^o)\otimes V(\bpi^e)$.  Analogous statements hold if $k$ is odd.
\end{prop}
\begin{pf} The first assertion of the proposition is immediate from Theorem \ref{hl}. However, it is easy to give a proof using Proposition \ref{cbraidthm} and we include it for completeness  since it is crucial for this paper.

 Define integers $r_j$, $1\le j\le k$, by\begin{gather*} r_1= 0, \ \ r_2=i_2-i_1+2,\\ r_{2s+1}= -i_1+2(i_2-i_3+\cdots -i_{2s-1}+i_{2s})-i_{2s+1},\ \ s\ge 1,\\
r_{2s+2}= -i_1+2(i_2-i_3+\cdots +i_{2s}-i_{2s+1})+i_{2s+2}+2,\ \ s\ge 1,\end{gather*} and note that $a_j=q^{r_j+m}$, for $1\le j\le k$.

We prove that $V(\bpi^o)$ is irreducible. Using Proposition \ref{cbraidthm}, it suffices to prove that for all $s>j$ and $p\in\bz_+$ with $p+1>i_{2s+1}$,   $$r_{2s+1}-r_{2j+1}\ne \pm(2p+2-i_{2j+1}-i_{2s+1}),$$ or equivalently that
 $$ -i_{2j+1}+2(i_{2j+2}-i_{2j+3}+\cdots -i_{2s-1}+i_{2s})-i_{2s+1}\ne \pm(2p+2-i_{2j+1}-i_{2s+1}).$$
This amounts to proving that
$$p+1\ne i_{2s}-i_{2s-1}+\cdots -i_{2j+3}+i_{2j+2} $$ and $$
p+1\ne i_{2s+1} - i_{2s} + \cdots+i_{2j+3}- i_{2j+2}+ i_{2j+1}.$$ If equality were to hold, then in the first case we would  get $p+1<i_{2s}<i_{2s+1}$ and in the second case we get $p+1<i_{2s+1}$ which contradicts our assumptions on $p$. The proof of the  irreducibility of $V(\bpi^o)$ is complete.  A similar argument proves the result for $V(\bpi^e)$.

 To prove that $V(\bpi^o)\otimes V(\bpi^e)$ contain $V(\bpi)$ as its unique submodule, we again use Proposition \ref{cbraidthm} and note that it is enough to   check that for all $s$ and $j$, we have
\begin{equation}\label{notin} r_{2j-1}-r_{2s}\notin\{ 2+ 2p-i_{2s}-i_{2j-1}:  \max\{i_{2j-1}, i_{2s}\}<  p +1\le \min\{i_{2j-1}+ i_{2s}, n+1\}\}.\end{equation}
For clarity, we  prove this by breaking   up the checking  into several cases.
If  $s\ge j\ge 1 $ and $i_{2s}+i_{2j-1}\le n+1$,
we have $$ r_{2s}-r_{2j-1}=i_{2s}+ i_{2j-1}+2-2(i_{2j-1}-i_{2j}+\cdots-i_{2s-2}+i_{2s-1}),$$ i.e.,
\begin{gather*}r_{2j-1}-r_{2s}= -i_{2s}-i_{2j-1}+2(-1+(i_{2j-1}-i_{2j})+\cdots +(i_{2s-1}-i_{2s}) +i_{2s}).
\end{gather*}
Since $(-1+(i_{2j-1}-i_{2j})+\cdots +(i_{2s-1}-i_{2s}) +i_{2s})<i_{2s}$, we see that equation  \eqref{notin} is satisfied. On the other hand, if $j\ge s\ge 1 $ and $i_{2s}+i_{2j-1}\le n+1$, we have \begin{gather*} r_{2s}-r_{2j-1} = i_{2s}+i_{2j-1}- 2(-1+(i_{2s}-i_{2s+1})+\cdots +(i_{2j-2}-i_{2j-1}) +i_{2j-1}).
\end{gather*}Since  the expression in parentheses is less than  than $i_{2j-1}$  we see that equation \eqref{notin} is again satisfied.
The other two cases are similar and we omit the details.\end{pf}

\subsection{Existence of $\varphi_2$}\label{missing}\hfill

 Using Proposition \ref{cyclic}  and the discussion in Section \ref{outline} we see that there exists a map $\varphi_2: L(\bpi)\to L(\bpi^o)\otimes L(\bpi^e)$. Moreover, since $\wt\bpi=\wt\bpi^o+\wt\bpi^e$, and $$L(\bpi)_{\wt\bpi}= \bc v_{\overline{\bpi}},\ \  L(\bpi^o)_{\wt\bpi^o}= \bc v_{\overline{\bpi^o}},\ \  L(\bpi^e)_{\wt\bpi^e}= \bc v_{\overline{\bpi^e}},$$ we see that $$\varphi_2( v_{\overline{\bpi}})= v_{\overline{\bpi^o}}\otimes v_{\overline{\bpi^e}}.$$ Moreover, $V(\bpi^o)$ and $V(\bpi^e)$ satisfy the conditions of Theorem \ref{primefund} and hence we have  the defining relations of $L(\bpi^o)$ and $L(\bpi^e)$.  We have now established the second step of the proof of Theorem \ref{first}, modulo the identification of $L(\bpi^o)$ and $L(\bpi^e)$ with the level one Demazure module. This will be done in Section \ref{ss:completion}.

\section{Level two Demazure modules in the tensor product\\  of Level one Demazure modules}\label{oddtpeven} We establish the third step in the proof of Theorem \ref{first}.

\subsection{Extended and affine Weyl groups}\label{affine}\hfill

 We use freely the notation established in Section \ref{sln}.   Let $\lie n^\pm$ be the subalgebra spanned by the elements $x^\pm_{i,j}$, $1\le i\le j\le n+1$,  and set $\lie b=\lie h\oplus\lie n^+$. Set $h_{i,j}=[x^+_{i,j}, x^-_{i,j}]$ and note that $h_{i,i}=h_i$.
Define   elements $x_0^\pm$ and $h_0$ of $\widehat{\lie g}$  by $$x_0^+=x^\mp_{1,n}\otimes t^{\pm 1},\ \ h_0=c-h_{1,n},$$  and set $$\widehat{\lie h}=\lie h\oplus\bc c\oplus\bc d ,\ \ \ \   \widehat{\lie n}^+= \lie g\otimes t\bc[t]\oplus\lie n^+,\ \  \ \ \widehat{\lie b}=\widehat{\lie n}^+\oplus\widehat{\lie h}.$$
We shall  regard an element of $\lie h^*$ as an element of $\widehat{\lie h}^*$ by setting it to be zero on $c$ and $d$.
  Define elements $\Lambda_i\in\widehat{\lie h}^*$, $0\le i\le n$,  by $$\Lambda_0(\lie h\oplus\bc d)=0,\ \ \Lambda_0(c)=1,\ \ \Lambda_i=\omega_i+\Lambda_0,\ \ 1\le i\le n,$$  and let  $\delta\in\widehat{\lie h}^*$ be given  by $\delta(\lie h\oplus\bc c)=0$, $\delta(d)=1$.  Let $\widehat P$  be the $\bz$--span of $\{\delta\}$  and  $\{ \Lambda_i: 0\le i\le n\}$ and $\widehat P^+$ is the $\bz$--span of $\delta$ and the $\bz_+$--span of $\Lambda_i$, $0\le i\le n$.

Let $s_i$, $0\le i
\le n$  be the  simple reflection of  the affine Weyl group $\widehat W$; recall that it  acts  on $\widehat{\lie h}^*$ and $\widehat{\lie h}$ by $$s_i(\mu) = \mu - \mu(h_i)\alpha_i,\quad s_i(h) = h - \alpha_i(h)h_i,\qquad \ h\in \widehat{\lie h},\ \ \  \mu\in \widehat{\lie h}^*.$$
 Note that $$w\in\widehat{W}\implies w(c)=c\ \ {\rm and}\ \ w(\delta)= \delta.$$  The Weyl group $W$ of $\lie g$  is the subgroup of $\widehat{W}$ generated by the elements $s_i$, $1\le i\le n$, and is isomorphic to the symmetric group on $n$ letters.
  Let $w_0\in W$ be the unique element of maximal length.  The action of  $W$ on $\widehat{\lie h}^*$ preserves $P$ and $Q$  and we have an isomorphism of groups $$\widehat W\cong W\ltimes Q.$$ The extended affine Weyl group $\widetilde W $ is the semi--direct product $W\ltimes P$. The affine Weyl group is a normal subgroup of $\widetilde W$ and if $\cal T$ is the group of diagram automorphisms of $\widehat{\lie g}$, we have $$\widetilde W\cong \cal T\ltimes \widehat W.$$  Here $\cal T$ is just isomorphic to the cyclic group of order $n+1$. Since $\cal T$ preserves $\widehat P$ and $\widehat P^+$, we see that $\widetilde W$ preserves $\widehat P$.  The following formulae make explicit the action of $\mu\in P$ on $\widehat{\lie h^*}$:$$ t_\mu(\lambda)= \lambda-(\lambda,\mu)\delta,\ \  \lambda\in\lie{h}^*\oplus\bc\delta,\ \ \ t_{\mu}(\Lambda_0)=\Lambda_0+\mu-\frac12(\mu,\mu)\delta.$$

\subsection{Integrable highest weight modules } \hfill

Recall that a weight module $V$ for $\widehat{\lie h}$ is one where $\widehat{\lie h}$ acts diagonally. Let $\wt V\subset\widehat{\lie h}^*$  be the set of eigenvalues for this action and, given $\mu\in\wt V$, let $V_\mu$ be the corresponding eigenspace.
For $\Lambda\in\widehat{P}^+$, let $V(\Lambda)$ be the irreducible highest weight integrable $\widehat{\lie g}$--module, which is  generated by an element $v_\Lambda$ with defining relations $$\widehat{\lie n}^+ v_\Lambda=0, \ \  \ \ hv_\Lambda= \Lambda(h)v_\Lambda,\ \ \  \ (x_i^-)^{\Lambda(h_i)+1}v_\Lambda=0,$$ where $h\in\widehat{\lie h}$ and $0\le i\le n$. Note that $\wt V(\Lambda)\subset\Lambda-\widehat Q^+$. (Here $\widehat Q^+$ is the $\bz_+$--span of the elements $\alpha_i$, $0\le i\le n$ and $\delta$).  It is easily seen that  for all $r\in\mathbb Z$, we have an isomorphism of $\widetilde{\lie g}$--modules, \begin{equation}\label{dshift} V(\Lambda-r\delta)\cong V(\Lambda).\end{equation}   The following proposition is well--known (see \cite[Chapters 10 and 11]{Kac} for instance).
\begin{prop}\label{tp}
\begin{enumerit} \item[(i)] Let $\Lambda\in \widehat P^+$.  We have$$\dim V(\Lambda)_{\mu}=\dim V(\lambda)_{w\mu},\ \ {\rm{ for \ all}}\ \ w\in\widehat W,\ \ \mu\in\widehat{\lie h}^*.$$ In particular $\dim V(\Lambda)_{w\Lambda}=1$ for all $w\in\widehat W$.
\item[(ii)] Given $\Lambda',\Lambda''\in\widehat P^+$, we have $$V(\Lambda')\otimes V(\Lambda'')\cong\bigoplus_{\Lambda\in\widehat P^+} \dim\left(\Hom_{\widehat{\lie g}}\left(V(\Lambda), V(\Lambda')\otimes V(\Lambda'')\right) \right)V(\Lambda).$$Moreover, \begin{equation}\label{tp1} \dim\Hom_{\widehat{\lie g}}\left(V(\Lambda), V(\Lambda')\otimes V(\Lambda'')\right)=\begin{cases} 1,\ \ \Lambda=\Lambda'+\Lambda'',\\ 0,\ \ \Lambda\notin\Lambda'+\Lambda''-\widehat Q^+.\end{cases}\end{equation} \end{enumerit}
\hfill\qedsymbol
\end{prop}
\begin{cor}
Suppose that  $\Lambda',\Lambda''\in\widehat P^+$ and that   $\Lambda=\Lambda'+\Lambda'' $.  For all  $w\in\widehat W$, we have \begin{equation}\label{tp2}\left (V(\Lambda')\otimes V(\Lambda'')\right)_{w\Lambda} = V(\Lambda)_{w\Lambda},\end{equation} where we have identified $V(\Lambda)$ with its image in $V(\Lambda')\otimes V(\Lambda'')$.\end{cor}
\begin{pf} Since the right  hand side is of dimension one by part (i) of the proposition and clearly contained in the left  hand side, it suffices to prove that $$\dim\left (V(\Lambda')\otimes V(\Lambda'')\right)_{w\Lambda} =1.$$ If not, we get by using part (ii) of the proposition with the first case of equation \eqref{tp1} that there exists $\Lambda_1\in\widehat P^+$ with $\Lambda_1\ne \Lambda$, such that $$\dim\Hom_{\widehat{\lie g}}\left(V(\Lambda_1), V(\Lambda')\otimes V(\Lambda'')\right)\ne 0,\ \qquad  V(\Lambda_1)_{w\Lambda}\ne 0.$$  Using part (i) of the proposition, this means that  $V(\Lambda_1)_\Lambda\ne 0$. But this is impossible, since $\Lambda_1\subset\Lambda-\widehat Q^+$ and $\Lambda_1\ne \Lambda$ thus  proving the corollary.
\end{pf}

\subsection{Stable Demazure modules} \hfill

Given $\Lambda\in\widehat P^+$ and $w\tau\in\widetilde W$, where $w\in\widehat W$ and $\tau\in\cal T$, the Demazure module $V_{w}(\tau\Lambda)$ is  the $\widehat{\lie b}$--submodule of $V(\tau\Lambda)$ given by
$$V_w(\tau\Lambda)=\bu(\widehat{\lie b})v_{w\tau\Lambda}, \ \  0\ne v_{w\tau\Lambda}\in V(\tau\Lambda)_{w\tau\Lambda}.$$  The Demazure modules are necessarily finite--dimensional since $\wt V(\lambda)\subset \Lambda-Q^+$. We say that $V_w(\tau\Lambda)$ is a level $\ell$--Demazure module if $\Lambda(c)=\ell$.  The following is immediate from Corollary \ref{tp}.

\begin{lem}\label{top} Let $w\tau\in\widetilde W$ and $\Lambda',\Lambda''\in\widehat P^+$. We have an isomorphism of $\widehat{\lie b}$--modules,  $$V_w(\tau(\Lambda'+\Lambda''))\cong \bu(\widehat{\lie b})(v_{w\tau\Lambda'}\otimes v_{w\tau\Lambda''})\subset V(\tau\Lambda')\otimes V(\tau\Lambda'').$$\hfill\qedsymbol
\end{lem}

{\em  In  this paper we are only interested in the  Demazure modules $V_w(\Lambda)$ satisfying the condition
   $w\Lambda(h_i)\le 0$, for all $1\le i\le n$.}   In this case, we have  $\lie n^- v_{w\Lambda}=0$ and $V_w(\Lambda)$ is a  module for the parabolic subalgebra  $\widehat{\lie b}\oplus\lie n^-$, i.e.,  $$V_w(\Lambda)=\bu(\widehat{\lie b}\oplus\lie n^-) v_{w\Lambda}=\bu(\lie g[t]) v_{w\Lambda} =\bu(\lie g[t]) v_{w_0^{-1}w\Lambda},$$  where the last equality follows from the fact that $V_w(\Lambda)$ is finite--dimensional $\lie g$--module. Writing  $w\Lambda=w_0\lambda+\Lambda(c)\Lambda_0+r\delta$, for a unique $\lambda\in P^+$ and $r\in\mathbb Z$, we see from equation \eqref{dshift} that $$V_w(\Lambda) \cong_{\lie g[t]} V_w(\Lambda-r\delta).$$ Hence, we denote these modules as $\tau_r^* D(\ell,\lambda)$, where $\ell=\Lambda(c)$, $r\in\bz$. Notice that the action of $d$ on these modules defines a $\bz$--grading on them which is compatible with the grading on $\lie g[t]$, i.e., the modules $\tau_r^* D(\ell,\lambda)$ are graded $\lie g[t]$--modules; for a fixed $\ell$ and $\lambda$ these modules are just grade shifts,  and we set $\tau^*_0D(\ell,\lambda)= D(\ell,\lambda)$. The eigenspace  $D(\ell,\lambda)_\lambda$ for the  $\lie h$--action is one--dimensional and we shall frequently denote a non--zero element of this space by $v_{\ell,\lambda}$.

\subsection{The main result on  tensor products of level one Demazure modules}\hfill

 The main result of this section is,
\begin{thm}\label{main1}  Given $\mu\in P^+$, there exists  $\mu^o,\mu^e\in P^+$ with $\mu=\mu^o+\mu^e$  such that we have an injective map of  graded $\lie g[t]$--modules $$D(2, \mu)\hookrightarrow D(1,\mu^o)\otimes D(1,\mu^e),\ \  v_{2,\mu}\to v_{1,\mu^o}\otimes v_{1, \mu^e}.$$
\end{thm}

\subsection{The key Proposition}\hfill

  Recall that   $$P^+(1)=\{\lambda\in P^+: \lambda(h_i)\le 1\ \ {\rm{for\ all}}\ \  1\le i\le n\}.$$ Given $\lambda=\sum_{j=1}^k\omega_{i_j}\in P^+(1)$, with $1\le i_1<i_2<\cdots<i_k\le n$, define $\lambda^o,\lambda^e\in P^+(1)$ by:
$$ \lambda^o=\begin{cases}\omega_{i_1}+\omega_{i_3}+\cdots+\omega_{i_k},\ \ \ k \ {\rm odd}, \\ \omega_{i_1}+\omega_{i_3}+\cdots+\omega_{i_{k-1}},\ \ \ k \ {\rm even},\end{cases}\ \qquad
\lambda^e=\lambda-\lambda^o.$$
We shall prove,
\begin{prop}\label{thm1} Given  $\lambda\in P^+(1)$ and $\nu\in P^+$ there exists  $w\in \widetilde W$ such that \begin{gather*} w(\nu+\lambda^o+\Lambda_0)\in \widehat P^+, \ \  w(\nu+\lambda^e+\Lambda_0)\in \widehat P^+.\end{gather*}
\end{prop}
\subsection{Proof of Theorem \ref{main1} and the third step of the proof of Theorem \ref{first}}\hfill

 Assuming Proposition \ref{thm1} the proof of Theorem \ref{main1} is completed as follows.  Write $\mu=2\nu+\lambda$, where $\nu\in P^+$ and $\lambda\in P^+(1)$, and set $\mu^o=\nu+\lambda^o$ and $\mu^e=\nu+\lambda^e$.  Choose  $w\in\tilde W$ as  in Proposition \ref{thm1}  and take  $$\Lambda =w(\mu+2\Lambda_0),\ \ \ \Lambda^o= w(\nu+\lambda^o+\Lambda_0), \ \  \Lambda^e=w(\nu+\lambda^e+\Lambda_0).$$ Then $\Lambda^o,\Lambda^e\in\widehat P^+$ and $\Lambda=\Lambda^o+\lambda^e\in\widehat P^+$ and  $$D(2,\mu) =V_{w_0w^{-1}}(\Lambda),\ \ D(1,\mu^o)= V_{w_0w^{-1}}(\Lambda^o),\ \  D(1,\mu^e)= V_{w_0w^{-1}}(\Lambda^e).$$ Theorem \ref{main1} is now immediate from Lemma \ref{top}.

The third step of the proof of  Theorem \ref{first} now follows.  Given $\bpi\in\cal P^+(1)$,  observe that $\wt\bpi^o=(\wt\bpi)^o$ and $\wt\bpi^e=(\wt\bpi)^e$ and hence using Theorem \ref{main1} we have a map of $\lie g[t]$--modules  $D(2,\wt\bpi)\hookrightarrow D(1,\wt\bpi^o)\otimes D(1,\wt\bpi^e)$.

\subsection{Proof of Proposition \ref{thm1}; reduction to $\nu=0$.} \hfill

In the  rest of this section we prove Proposition \ref{thm1}.  The first step is to show that, for a fixed $\lambda$, it suffices to prove the result when $\nu=0$. Thus, suppose that we have chosen $w\in\widetilde W$ such that $w(\lambda^o+\Lambda_0)$ and $w(\lambda^e+\Lambda_0)$ are in $\widehat P^+$.   Since $(\lambda^o+ \Lambda_0)(c)=1=(\lambda^e+\Lambda_0)(c)$, we may write $$w(\lambda^o+\Lambda_0)=\Lambda_i+p^o\delta,\ \  w(\lambda^e+\Lambda_0)=\Lambda_j+p^e\delta,$$ for some $p^o,p^e\in\mathbb Z$ and $0\le i,  j\le n$.   Using the formulae in Section \ref{affine}, and the fact that $\Lambda_p=\omega_p+\Lambda_0$ for $0\le p\le n$, where $\omega_0=0$, we get $$t_{-w\nu}w(\lambda^o+\Lambda_0+\nu)=\left( \Lambda_i +(p^o+\frac 12(\nu,\nu)+(\omega_i, w\nu))\delta\right)\in\widehat P^+,$$ $$t_{-w\nu}w(\lambda^e+\Lambda_0+\nu)=\left( \Lambda_j +(p^e+\frac 12(\nu,\nu)+(\omega_j, w\nu))\delta\right)\in\widehat P^+,$$ and the claim is established.

\subsection{Proof of Proposition \ref{thm1};  the case $\nu=0$.} \hfill

 Consider  the partial order on $P^+(1)$ given by $\mu\le \nu$ iff $\nu-\mu\in Q^+$.  The minimal elements of this order are $0$ and $\omega_i$, $1\le i\le n$. If $\lambda=0$ the result is clear, we just take $\Lambda^o=\Lambda^e=\Lambda_0$ and $w=\id$. If  $\lambda=\omega_i$, we take $\lambda^o=\omega_i$ and $\lambda^e=0$. Since  $\omega_i+\Lambda_0\in\widehat P^+$ we take $\Lambda^o=\omega_i+\Lambda_0$, $\Lambda^e=\Lambda_0$ and $w=\id$. If $\lambda=\omega_i+\omega_j$ with $i<j$  then we again take $w=\id$ since $\omega_p+\Lambda_0=\Lambda_p$ for all $0\le p\le n$.

For the inductive step,  let $\lambda\in P^+(1)$  with $\lambda=\sum_{j=1}^k\omega_{i_j}$, $i_1<\cdots <i_k$ and $k>2$. Suppose that we have proved the result for all elements $\mu\in P^+(1)$ with $\mu<\lambda$. To prove the result for $\lambda$ it  clearly suffices to show that there exists $w\in\widetilde W$ and $\mu\in P^+$ with $\mu<\lambda$, and $p^o,p^e\in\mathbb Z$, such that  \begin{equation}\label{e:redto} w(\lambda^o+\Lambda_0)=\mu^o+\Lambda_0+p^o\delta,\ \ {\rm and}\ \ w(\lambda^e+\Lambda_0)=\mu^e+\Lambda_0+p^e\delta.\end{equation}  This is done as follows: take $$w=\begin{cases}
s_{i_3}s_{i_3+1}\cdots s_ns_0,\ \ k=3, \ \  i_1=1,\\  s_{i_3}s_{i_3+1}\cdots s_ns_{i_{k-2}-1}s_{i_{k-2}-2}\cdots s_1s_0,\ \ k >3,\ \ {\rm{or}}\ \ k=3,\ \  i_1>1,\end{cases}$$
and $$\mu=\begin{cases} \omega_{i_1-1}+\omega_{i_2}+ \omega_{i_3+1},\ \ k=3,\\
\omega_{i_1-1}+\omega_{i_2-1}+\omega_{i_3}+\cdots +\omega_{i_{k-2}}+\omega_{i_{k-1}+1}+\omega_{i_k+1},\
 \ \ k>3. \end{cases}$$ Note that $\mu<\lambda$ since,  $$ \lambda-\mu=\begin{cases} \alpha_{i_1}+\cdots+\alpha_{i_3},\ \ k=3,\\ \alpha_{i_1}+2\alpha_{i_2}+\cdots+2\alpha_{i_{k-1}}+\alpha_{i_k}),\ \ k>3.\end{cases}$$
We now establish that   equation \eqref{e:redto} is satisfied. For this, it is most convenient to deal with the cases $k= 3,4$ separately.
 If  $k=3$ and $i_1>1$ or if $k=4$, a simple calculation gives   \begin{gather*} w(\lambda^o+\Lambda_0)=w(\Lambda_{i_1}+\Lambda_{i_3}-\Lambda_0) = \Lambda_{i_1}+\Lambda_{i_3}-\Lambda_0+\sum_{j=0}^{i_1-1}\alpha_j+\sum_{j=i_3}^n\alpha_j=\mu^o+\Lambda_0+\delta.\end{gather*} Moreover, if  $k=3$, we have  $$w(\lambda^e+\Lambda_0)=w(\Lambda_{i_2})= \Lambda_{i_2}=
\mu^e+\Lambda_0,$$ while if $
k=4$, we have $$ w(\lambda^e+\Lambda_0)=w(\Lambda_{i_2}+\Lambda_{i_4}-\Lambda_0)=  \Lambda_{i_2}+\Lambda_{i_4}-\Lambda_0+\sum_{j=0}^{i_2-1}\alpha_j+\sum_{j=i_4+1}^n\alpha_j=\mu^e+\Lambda_0+\delta.$$
The case $k=3$ and $i_1=1$ is identical and we omit the details.
  For $k\ge 5$,  we write  \begin{gather*}\label{wdelta}w(\lambda^o+\Lambda_0)= w(\lambda^o+\Lambda_0)(d)\delta+ \sum_{j=0}^n (w(\lambda^o+\Lambda_0), \alpha_j)\Lambda_j
\end{gather*} and similarly for $\lambda^e+\Lambda_0$. We have  $$w(\lambda^o+\Lambda_0)(d) =(\lambda^o+\Lambda_0)(w^{-1}d)=(\lambda^o+\Lambda_0)(d-h_0) =\lambda^o(h_\theta)-1.$$ To prove that  $(w(\lambda^o+\Lambda_0),\alpha_j)=(\mu^o+\Lambda_0,\alpha_j)$ it is enough to prove that  $$(\lambda^o+\Lambda_0)(w^{-1}h_j)=(\mu^o+\Lambda_0)(h_j)$$ and this is done by using the following easily established formulae:
 \begin{gather*}w^{-1}(\alpha_0)=\alpha_0+\alpha_1+\alpha_n,\qquad \ w^{-1}(\alpha_j)=\alpha_j,\ \ i_3<j<i_{k-2},\\
w^{-1}(\alpha_j)= \alpha_{j+1},\ \  0<j<i_3-1,\qquad\ w^{-1}(\alpha_j)= \alpha_{j-1},\ \ j>i_{k-2}+1,\\
w^{-1}(\alpha_{i_3-1})=\alpha_{i_3}+\cdots +\alpha_n+\alpha_0,\qquad  w^{-1}(\alpha_{i_{k-2}+1})= \alpha_{i_{k-2}}+\cdots+\alpha_1+\alpha_0,\\
w^{-1}(\alpha_{i_{k-2}})= -(\alpha_0+\alpha_1+\cdots +\alpha_{i_{k-2}-1}),\end{gather*}
and
\begin{gather*}w^{-1}(\alpha_{i_3})=\begin{cases}-(\alpha_{i_3+1}+\cdots+\alpha_n+\alpha_0),\ \ k>5,\\  -(\alpha_1+\cdots+\alpha_n+2\alpha_0),\ \ k=5.\end{cases}\end{gather*} The case of $w(\lambda^e+\Lambda_0)$ is identical and the proof of Proposition \ref{thm1} is complete.

\section{A presentation of level two $\lie g$--stable  Demazure modules.} \label{effpres}

 An infinite   set of    defining relations  of the Demazure module $V_w(\Lambda)$ was  given in \cite{Joseph, Mathieu} for all $w\in\widetilde W$ and $\Lambda\in\widehat P^+$. In the case when $w\Lambda(h_i)\le 0$, for all $1\le i\le n$, these relations can be used (see \cite{FoL}, \cite{Naoi})  to give (a still infinite set) of defining relations  of $V_w(\Lambda)$  (or equivalently of $\tau_r^*D(\ell,\lambda)$) as  $\lie g\otimes \bc[t]$--modules.  In the case of the level one Demazure modules, it was shown in \cite{CL} for $\lie{sl}_{n+1}$,  that the relations could be reduced to a finite number of relations. This was later extended for arbitrary levels in \cite{CV}. In this section we show that, in the case of $A_n$ and $\ell=2$, one can further whittle down the  set of defining relations and as a consequence, we  prove Proposition \ref{second}. Along the way, we shall see that if $\bpi\in\cal P^+(1)$, then $L(\bpi^o)$ and $L(\bpi^e)$ are level one Demazure modules, which establishes the missing piece (see Section \ref{missing})  of the second step of Theorem \ref{first}.

\subsection{A refined presentation of Demazure modules}\label{ss:completion} \hfill

 We recall the presentation of $D(\ell,\mu)$, $\mu\in P^+$, given in \cite[Theorem 2]{CV}.
\begin{prop}\label{demrelref} For $\ell\in\mathbb N$,  $\mu\in P^+$, the $\lie g[t]$--module  $D(\ell,\mu)$ is generated by an element $w_\mu$ satisfying the following defining relations:\begin{gather*}(x_i^+\otimes 1)w_\mu=0,\ \ \ (h_i\otimes t^r)w_\mu=\mu(h_i)\delta_{r,0}w_\mu,\ \ \ \ (x_i^-\otimes 1)^{\mu(h_i)+1}w_\mu=0,\\  (x^-_{i,j}\otimes t^{s_{i,j}}) w_\mu=0,\ \  1\le i\le j\le n, \\
(x^-_{i,j}\otimes t^{s_{i,j}-1})^{m_{i,j}+1}w_\mu=0, \end{gather*}
where $r\in\bz_+$,  $1\le i\le j\le n$ and $s_{i,j},m_{i,j}\in\bz_+$ are uniquely defined by requiring $\mu(h_{i,j})= (s_{i,j}-1)\ell+m_{i,j}$ with  $0<m_{i,j}\leq\ell$. Moreover, if $m_{i,j}=\ell$, the final  relation is a consequence of the preceding  relations.
\hfill\qedsymbol
\end{prop}
\begin{cor} For all $\ell\in\bn$ and $\mu\in P^+$, the assignment $w_\mu\to w_\mu$ defines a canonical surjective map $D(\ell,\mu)\to D(\ell+1,\mu)\to0$ of $\lie g[t]$--modules. \hfill\qedsymbol
\end{cor}
\begin{rem}   It is clear from Lemma \ref{lpi} that, for all $\bpi\in\cal P^+_\bz$, the module $L(\bpi)$ is a quotient of $D(1,\wt\bpi)$.   If $\bpi\in\cal P^+(1)$, it follows from Theorem \ref{primefund} and  Proposition \ref{cyclic} that $$L(\bpi^o)\cong D(1,\wt\bpi^o), \ \ L(\bpi^e)\cong D(1,\wt\bpi^e),$$ which establishes  missing piece (see Section \ref{missing}) of the proof of the second step of Theorem \ref{first}.
\end{rem}

{\subsection{Evaluation modules}\hfill

   Let $\ev_0:\lie g[t]\to\lie g$ be the map of Lie algebras given by setting $t=0$,i.e., $a\otimes f\to  f(0)a$, for all $a\in\lie g$ and $f\in\bc[t]$. Given a $\lie g$--module $V$, let  $\ev_0^*V$  denote the corresponding  $\lie g[t]$--module. The following is straightforward.
\begin{lem}\label{elem} For all $\ell\in\bz_+$ and $\mu\in P^+$, we have $\wt D(\ell,\mu)\subset\mu-Q^+$. Moreover, $$\dim\Hom_{\lie g}(V(\mu), D(\ell,\mu))=\dim\Hom_{\lie g[t]}(D(\ell,\mu),\ev_0^* V(\mu))= 1, $$  and$$D(\ell,\lambda)\cong_{\lie g[t]} \ev_0^*V(\lambda),\ \ {\text{if}}\ \ \lambda(h_{1,n})\le \ell. $$\hfill\qedsymbol \end{lem}}

\subsection{The $\lie{sl}_2$ case}\hfill

In the case of $n=1$, i.e., when $\lie g$ is of type $\lie{sl}_2$, we   identify $P^+$ with $\bz_+$ freely and let $x^\pm $  be the root vectors $x_1^\pm$.  Given a partition, $\xi=(\xi_1\ge\cdots\geq\xi_m\ge 0)$, $r\in\bn$, define a $\lie{sl}_2[t]$--module $V(\xi)$ as follows: it is the cyclic module generated by an element $v_\xi$ with defining relations:
\begin{equation}\label{vxirel1} (x^+\otimes t)v_\xi =0, \ \ (h\otimes t^r)v_\xi=|\xi|\delta_{r,0}v_\xi,\ \ (x^-\otimes 1)^{|\xi|+1}v_\xi=0,\ \ |\xi|=\sum_{k\ge 1}\xi_k,\end{equation}  and \begin{equation}\label{vxirel2}(x^+\otimes t)^s(x^-\otimes 1)^{s+r}v_\xi=0,\end{equation} for all $s,r$ satisfying the condition, that there exists $k\in\bz_+$ with $s+r\ge 1+rk+\sum_{p\ge k+1}\xi_k.$

In this paper, we shall only be interested in the case when  the maximum value of a part is two, i.e., $\xi$ is of the form $2^a1^b$ for some $a,b\in\bz_+$.
We summarize for this special case, the results of \cite[Theorem 2, Theorem 5]{CV} that we shall need.
\begin{prop}\label{vxi} Let $\xi=2^a1^b$, $a,b\in\bz_+$  be a partition.
\begin{enumerit}
\item[(i)] We have  $V(1^b)\cong D(1, b)$ as $\lie{sl}_2[t]$--modules.
\item[(ii)]  If $b\in\{0,1\}$, then $V(2^a1^b)\cong D(2,2a+b)$ as $\lie{sl}_2[t]$--modules.
\item[(iii)] If $b\ge 2$, then we have a short exact sequence of graded $\lie g[t]$--modules,\begin{equation}\label{ses} 0\to V(2^a 1^{b-2})\stackrel{\iota}{\to }V(2^a 1^b)\stackrel{\pi}{\to} V(2^{a+1} 1^{b-2})\to 0, \end{equation} such that $$\iota(v_{2^a1^{b-2}})= (x\otimes t^{a+b-1})v_{2^a1^b},\qquad  \pi(v_{2^a1^b})= v_{2^{a+1}1^{b-2}}.$$
\end{enumerit} \hfill\qedsymbol
\end{prop}
\begin{cor} The module  $V(2^a1^b)$ is generated by the element $v_{2^a1^b}$ with the relations given in \eqref{vxirel1} and the single additional relation: $(x^-\otimes t^{a+b})v_{2^a1^b}=0$.
\end{cor}
\begin{pf} If $a=0$, then the corollary follows from Proposition \ref{demrelref} and part (i) of Proposition \ref{vxi}.  A straightforward induction on $a$ together with the short exact sequence in \eqref{ses} establishes the corollary.
\end{pf}

\subsection{A further refinement of the presentation of level two Demazure modules} \

\begin{prop} \label{sl2}\label{leveltwo} Assume that $\lie g$ is of type $\lie{sl}_{n+1}$ and let   $\mu\in P^+$. For $1\le i\le j\le n$, write $\mu(h_{i,j})=2(s_{i,j}-1)+m_{i,j}$ with  $m_{i,j}\in\{1,2\}$.  The relation $$
(x^-_{i,j}\otimes t^{s_{i,j}-1})^{m_{i,j}+1}w_\mu=0$$ in $D(2,\mu)$ is redundant.\hfill\qedsymbol

\end{prop}
\begin{pf} It  follows from Proposition \ref{demrelref} that it suffices to prove the proposition when $m_{i,j}=1$.   If  $n=1$, then the result follows from Proposition \ref{vxi}(ii) and Corollary \ref{vxi}.  Otherwise, set $\alpha=\alpha_{i,j}$ and consider the subalgebra $\lie s_\alpha[t]$ of $\lie g[t]$  spanned by the elements $x^\pm_\alpha\otimes\bc[t]$.  Clearly $\lie s_\alpha[t]$ is isomorphic to $\lie {sl}_2[t]$ with the element $x^\pm_\alpha$ mapping to $x^\pm$. Moreover, if we denote by$D_{\lie{s}_2[t]}(2,\mu(h_\alpha))$ the level two Demazure module for $\lie{sl}_2[t]$, then we have a non--zero map  map of $\lie{sl}_2[t]$--modules, $$D_{\lie{sl}_2[t]}(2,\mu(h_\alpha))\to \bu(\lie{ s}_\alpha[t])w_\mu\subset D(2,\mu).$$ The proposition is now immediate from the $n=1$ case.
\end{pf}

\subsection{Proof of Proposition \ref{second}} We prove Proposition \ref{second}.   Let $\lambda=\omega_{i_1}+\cdots+\omega_{i_k}\in P^+(1)$,where $1\leq i_1 < \cdots < i_k\leq n$ and $\nu\in P^+$.  For $\alpha\in R^+$, write $2\nu+\lambda(h_\alpha)= 2(s_\alpha-1)+m_\alpha$ with $s_\alpha\ge 1$ and $m_\alpha\in\{1,2\}$.  Note that $$ \nu(h_i)+\lambda(h_i)= s_i,\ 1\le i\le n,\ \ \ \ \ \nu(h_{i_p,i_{p+1}})+1=s_{i_p,i_{p+1}},\ \ \ 1\le p\le k-1.$$
It is clear from the definition of $M(\nu,\lambda)$ (see Section \ref{dnulambda}) and Proposition \ref{demrelref} that there exists a sujrective  morphism $$M(\nu,\lambda)\to D(2,2\nu+\lambda)\to 0$$ of $\lie g[t]$--modules.
The proposition follows if we  prove that the preceding map is an isomorphism.
 By Corollary  \ref{leveltwo}, it suffices to show that \begin{equation}\label{toprove} (x^-_{i_p,i_{p+1}}\otimes t^{s_{i_p, i_{p+1}}})v_{\nu,\lambda}=0,\ \ 1\le p\le k-1 \implies (x^-_{i,j}\otimes t^{s_{i,j}})v_{\nu,\lambda}=0,\ \ \  1\le i\le j\le n.\end{equation}
 We proceed by induction on $j-i$, with  induction beginning at $i=j$ by  hypothesis (see \eqref{def2}). Assume that we have proved the result for all $\alpha_i+\cdots+\alpha_j$ with $j-i<s$ and consider the case when $j=i+s$.  Suppose first that there does not exist $1\le p\le k$ such that $i\le i_p\le i+s$. In this case, we have $\mu(h_{i,j})=2\nu(h_{i,j})=2\nu(h_i)+ 2\nu(h_{i+1, j})$.   The induction hypothesis implies that $$(x^-_{i}\otimes t^{\nu(h_i)}) v_{\nu,\lambda}=0=(x^-_{{i+1,j}}\otimes t^{\nu(h_{i+1,j})}) v_{\nu,\lambda}.$$ Since $[x^-_{\alpha_i}, x^-_{\alpha_{i+1,j}} ]=x^-_{\alpha_{i,j}}$,  we get \eqref{toprove} in this case.
We consider the other case, when we can choose   $1\le p\le k$ minimal and $1\le r\le k$ maximal so that $i\le i_p\le i_r\le j$. If $i<i_p$ we have again by the inductive hypothesis that  $$(x^-_{i}\otimes t^{\nu(h_i)}) v_{\nu,\lambda}=0=(x^-_{i+1,j}\otimes t^{(\nu+\lambda)(h_{i+1,j})}) v_{\nu,\lambda}$$  and the inductive step is completed as before. If $j>i_r$ then the proof is similar where we write $\alpha_{i,j}=\alpha_{i,j-1}+\alpha_{j}$. Finally, suppose that $i=i_p$ and $ j=i_r$. If $r=p+1$ then the inductive step is the hypothesis in equation \eqref{def3}. If  $r\ge p+2$, then we write $\alpha_{i,j}=\alpha_{i_p,i_{p+1}}+\alpha_{i_{p+1}+1, i_r}$. This time the induction hypothesis gives,
$$(x^-_{{i_p,i_{p+1}}}\otimes t^{\nu(h_{i_p,i_{p+1}})+1})v_{\nu,\lambda}\ \ =\  0\ =\ \ (x^-_{{i_{p+1}+1,i_{r}}}\otimes t^{(\nu+\lambda)(h_{i_{p+1}+1,i_{r}})})v_{\nu,\lambda},$$ and the inductive step is completed as  before.

\section{A characterization of $\lie g$--stable level two Demazure modules}\label{fusionproducts}
We establish the final step of the proof of Theorem \ref{first}.

\subsection{}  We shall prove the following result in the rest of the section. \begin{prop} \label{gengradedlimit}  Let $\mu\in P^+$ and let  $V$ be a  (not necessarily graded) $\lie g[t]$--module which is isomorphic to $D(2,\mu)$ as a $\lie g$--module. Assume that
$\varphi: D(1,\mu)\to V\to 0$ is a surjective map of $\lie g[t]$--modules. Then $V$ is isomorphic to $D(2,\mu)$ as $\lie g[t]$--modules.
\end{prop}  Let $\nu\in P^+$ and $\lambda=\omega_{i_1}+\cdots +\omega_{i_k}\in P^+(1)$, $1\le i_1<\cdots <i_k\le n$, be the unique elements such that $\mu=2\nu+\lambda$.   Setting  $\varphi(w_\mu)=v_\mu\in V $, we see that $v_\mu$ satisfies the relations in  \eqref{def1}.
Since $\dim V=\dim D(2,\mu)$ by hypothesis, it follows from Proposition \ref{second} that it suffices to prove that the element $v_\mu $, satisfies \eqref{def2} and \eqref{def3}.
\subsection{}
\begin{lem}\label{wtspacedim1} For $1\le i\le n$, with $(\nu+\mu)(h_i)>0$, we have $\dim V_{\mu-\alpha_i}\le (\nu+\lambda)(h_i).$\end{lem}

\begin{pf} Since $$V\cong_{\lie g}  D(2,\mu)\implies V_{\mu-\alpha_i}\cong D(2,\mu)_{\mu-\alpha_i},$$ it suffices to prove that $\dim  D(2,\mu)_{\mu-\alpha_i}\le (\nu+\lambda)(h_i).$  A simple application of the Poincare--Birkhoff-Witt (PBW) theorem, along with the fact that $(x^-_i\otimes t^{(\nu+\lambda)(h_i)})w_{\mu}=0$ in $D(2,\mu)$, shows that $D(2,\mu)_{\mu-\alpha_i}$ is spanned by the elements $\{(x_i^-\otimes t^r)w_{\mu}: 0\le r<(\nu+\lambda)(h_i)\}$ and hence proves the Lemma. \iffalse The non--zero elements of this set all have different grades and it follows that they are linearly independent elements of $D(2,\mu)_{\nu+\lambda-\alpha_i}$ which proves
$$\dim D(2,\mu)_{\mu-\alpha_i}\leq(\nu+\lambda)(h_i),$$ as required.\fi
\end{pf}
\subsection{} \label{reldef1}We now prove that $v_\mu$ satisfies equation \eqref{def2}. % We proceed by contradiction.  
Proposition \ref{demrelref} implies that the element  $w_\mu\in D(1,\mu)$  satisfies, $(x_i^-\otimes t^{\mu(h_i)})w_\mu=0$, $1\le i\le n$,  and  hence  we also have \begin{equation}\label{existsm}(x_i^-\otimes t^{\mu(h_i)})v_\mu=0, \ \ \ 1\le i\le n.\end{equation}  If  $(\nu+\lambda)(h_i)\in\{0,1\}$, then $\mu(h_i)=(\nu+\lambda)(h_i)$  and  there is nothing to prove. 

 Suppose that $\mu(h_i)\ge 2 $ or equivalently $\nu(h_i)>0$.  Using the second relation in  \eqref{def1} we get,
   $$(h_i\otimes t)(x_i^-\otimes t^s) v_\mu= -2 (x_i^-\otimes t^{s+1}) v_\mu,$$  and hence, $$(x_i^-\otimes t^s) v_\mu=0\implies
(x_i^-\otimes t^{s+1}) v_\mu=0.$$  We now proceed by contradiction to prove that $v_\mu$ satisfies \eqref{def2}. Thus, if  $(x^-_{i}\otimes t^{(\nu+\lambda)(h_i)}) v_\mu\ne 0$ then
 the  elements of the set  $\{ (x^-_{i}\otimes t^{s})v_\mu: 0\le s\le  (\nu+\lambda)(h_i)\}$   are all non--zero  and  by Lemma \ref{wtspacedim1} must be  linearly dependent.  It follows that  there exists $0\le m< (\nu+\lambda)(h_i)$ such that we have a non--trivial linear combination, 
 \begin{equation}\label{lincomb} \sum_{s=m}^{(\nu+\lambda)(h_i)}z_s(x_i^-\otimes t^s) v_\mu =0,\ \ \   z_s\in\bc, \ \  m\le s\le (\nu+\lambda)(h_i), \ \   z_m\ne 0.\end{equation}  
Since $\mu(h_i)-m-1\ge \nu(h_i)>0$, we apply $(h_i\otimes t^{\mu(h_i)-m-1})$ to the preceding equation  and use \eqref{def1},\eqref{existsm},  along with the fact that  $z_m\ne 0$ to get 
$$(x^-_{i}\otimes t^{\mu(h_i)-1}) v_\mu=0.$$  If  $\mu(h_i)\in\{2,3\}$, then we have $\nu(h_i)=1$ and $(\nu+\lambda)(h_i)=\mu(h_i)-1$ and we have the desired contradiction.   Otherwise, we have $\mu(h_i)\ge 4$, i.e., $\nu(h_i)\ge 2$ and hence we get $\mu(h_i)-m-2\ge \nu(h_i)-1>0$.  Hence applying
$(h_i\otimes t^{\mu(h_i)-m-2})$ to the  expression in equation \eqref{lincomb},  we now get    $(x^-_{i}\otimes t^{\mu(h_i)-2}) v_\mu=0$.   If $\mu(h_i)\in\{4,5\}$,  we would have $(\nu+\lambda)(h_i)=\mu(h_i)-2$ which would again contradict our assumptions.
Further iterations of this argument  give a contradiction for all values of  $\mu(h_i)$.  Hence, we have  $(x^-_{i}\otimes t^{(\nu+\lambda)(h_i)})v_\mu= 0$, as required.

\subsection{} We need some additional notation to complete the proof that $v_\mu$ satisfies \eqref{def3}.  Recall that $\alpha_{i,j}=\alpha_{i}+\cdots +\alpha_j$.  Let $$\bu(\lie n^-[t])_{\alpha_{i_p,i_{p+1}}}=\{u\in \bu(\lie n^-[t]): [h, u] =-\alpha_{i_p,i_{p+1}}(h) u, \ \textrm{ for all}\ h\in \lie h\}.$$  Observe that, if $h\in\lie h$ and $r\in\bz_+$, then $$\bom\in \oplus_{p\in\bz_+}\bs[p]\implies [h\otimes t^r,\bom]\in  \oplus_{p\in\bz_+}\bs[p+r].$$  The elements 
$$(x_{\beta_s}^-\otimes t^{\ell_s})\cdots (x_{\beta_1}^-\otimes t^{\ell_1}),\ \  s\ge 1 ,
 \ell_s\in\bz_+, $$ and $\beta_j\in R^+$, $\ell_j\in\bz_+$ satisfying:
  $ \beta_1=\alpha_{i_p}+\alpha_{i_p+1}+\cdots +\alpha_m$, for some $m\ge i_p$,
$\beta_1+\cdots +\beta_r\in R^+$, for all $1\le r\le s$, $\beta_1+\cdots+\beta_s=\alpha_{i_p,i_{p+1}}$ are a basis for $ \bu(\lie n^-[t])_{\alpha_{i_p,i_{p+1}}}$. 
Let $S\subset \bu(\lie n^-[t])$ be the subset  of this basis with the additional restriction:  $ 0\le \ell_j\le N_{\beta_j}$, where \begin{equation}\label{Nj} N_{\beta_j}=\begin{cases}  \nu(h_{i_p,i_{p+1}})+1,\ \ s=1,\\ \nu(h_{\beta_j}), \ \  j\in\{1, s\},\ \  s>1, \\\nu(h_{\beta_j})-1,\ \ 1<j<s.  \end{cases}\end{equation}
Then $S$ is a finite,  linearly independent subset of this space and we  let $\bs$ be the span of $S$. Notice also that $\bs$ is graded and that \begin{equation}\label{gradeds}\bs[p]=0,\ \ p>\nu(h_{i_p,i_{p+1}})+1
\ \ \ {\rm and}\ \ \ \bs[\nu(h_{i_p,i_{p+1}})+1]= \bc(x^-_{i_p,i_{p+1}}\otimes t^{\nu(h_{i_p, i_{p+1}})+1}).\end{equation}

\subsection{}\label{vdef3}  We turn to the proof that $v_\mu$ satisfies \eqref{def3}.   Denote by $\tilde V$ the graded quotient of $D(1, \mu)$ by the $\lie g[t]$--submodule  generated by the set $$\{(x^-_{i}\otimes t^{(\nu+\lambda)(h_i)})w_\mu: 1\le i\le n\},$$ and let $\tilde v_\mu$ be the image of $w_\mu$ in $\tilde V$.  The results in section \ref{reldef1} show that $V$ is a $\lie g[t]$--quotient of $\tilde V$; the definition of $D(2,\mu)$ given in Proposition \ref{demrelref} shows that $D(2,\mu)$ is a graded $\lie g[t]$--quotient of $\tilde V$. Since $$(x_i^-\otimes t^{\nu(h_i)})\tilde v_\mu =0,\ \ i_p<i<i_{p+1},\ \qquad  (x_i^-\otimes t^{\nu(h_i)+1)})\tilde v_\mu =0,\  \  \ i\in \{i_p, i_{p+1}\},$$ taking repeated commutators gives \begin{gather} (x^-_{i_p,i_{p+1}}\otimes t^{\nu(h_{i_p,i_{p+1}})+2})\tilde v_\mu= 0 \ \ {\rm and} \nonumber \\ (x^-_{i_p, i}\otimes t^{\nu(h_{i_p,i})+1})\tilde v_\mu = (x^-_{j, i_{p+1}}\otimes t^{\nu(h_{j,i_{p+1}})+1})\tilde v_\mu= (x^-_{i,j}\otimes t^{\nu(h_{i,j})})\tilde v_\mu=0,\label{e:reptcom}\end{gather}
for all $i_p<i\leq j<i_{p+1}$.

The following is a straightforward consequence  of  the relations in equation \ref{e:reptcom} and the PBW theorem.
\begin{lem}\label{spanv} The  space $\tilde{V}_{\lambda-\alpha_{i_p,i_{p+1}}}$ is spanned by the elements $\{\bx \tilde v_\mu:\bx\in S\}$, and hence $$\dim \tilde  V_{\lambda-\alpha_{i_p,i_{p+1}}}\le |S|.$$
\hfill\qedsymbol
\end{lem}
\subsection{} We now prove,
\begin{lem}\label{dimless} We have  $\dim   V_{\lambda-\alpha_{i_p,i_{p+1}}}< |S|$, i.e., there exists a non--zero element $\bom\in\bs$ such that $\bom v_\mu=0$. \end{lem}
\begin{pf} Assume for a contradiction that  $\dim   V_{\lambda-\alpha_{i_p,i_{p+1}}}= |S|$. Then we have $$\dim \tilde  V_{\lambda-\alpha_{i_p,i_{p+1}}} =|S|= \dim D(2,\mu)_{\lambda-\alpha_{i_p,i_{p+1}}},$$ where the first equality follows because   $V$ is a quotient of $\tilde V$ and the  second equality  follows because $V\cong D(2,\mu)$ as $\lie g$--modules. Since $D(2,\mu)$ is a quotient of $\tilde V$, this means that using Lemma \ref{spanv},  the elements $\{\bx\tilde v_\mu: \bx\in S\}$ and hence the elements $\{\bx w_\mu:\bx\in S\}$ are linearly independent subsets of $\tilde V$ and $D(2,\mu)$ respectively. But the latter is impossible since $$(x^-_{i_p,i_{p+1}}\otimes t^{\nu(h_{i_p,i_{p+1}})+1})\in S,$$ and $ (x^-_{i_p,i_{p+1}}\otimes t^{\nu(h_{i_p,i_{p+1}})+1})w_\mu=0$ is a defining relation  in $D(2,\mu)$. The lemma is proved.
\end{pf}

\subsection{} 

\begin{prop} Let $0\ne \bom\in\bs$ be such that $\bom v_\mu=0$ and assume that the minimal graded component of $\bom$ is  $k$, for some $k<\nu(h_{i_p, i_{p+1}})+1$. There exists $0\ne \bom'\in\bs$ with $\bom' v_\mu=0$ whose minimal graded component is at least $k+1$. In particular, we have  $$\bs[\nu(h_{i_p,i_{p+1}})+1]v_\mu= \bc(x^-_{i_p,i_{p+1}}\otimes t^{\nu(h_{i_p, i_{p+1}})+1})v_\mu=0$$ and hence $v_\mu$ satisfies \eqref{def3}.

\end{prop}
\begin{pf}  For  $1\leq i \leq n$, let $h_{\omega_i}\in \lie h$ be the unique element such that $\alpha_j(h_{\omega_i})=\delta_{i,j}$, for all $1\leq j \leq n$.  Write $\bom=\sum_{\bx\in S}z_\bx\bx$, $z_\bx\in\bc$, and note that $[h\otimes t^r,\bom]v_\mu=0$, for all $h\in\lie h$, $r\in\bz_+$.   Moreover,  we have
$$[h_{\omega_1}\otimes t, \bx]= (x_{\beta_s}^-\otimes t^{\ell_s})\cdots (x_{\beta_1}^-\otimes t^{\ell_1+1}),\quad \bx\in S.$$ Suppose that the set
$$S_1(\bom)= \{\bx= (x_{\beta_s}^-\otimes t^{\ell_s})\cdots (x_{\beta_1}^-\otimes t^{\ell_1})\in S: z_\bx\ne 0,\ \ \ell_1<N_{\beta_1}\}\ne\emptyset.$$Then,
the elements $ \{[h_{\omega_1}\otimes t, \bx]:\bx \in S_1(\bom)\}$ are all distinct elements of $S$
and, by \eqref{e:reptcom} and the definition of $N_{\beta_1}$, we have
 $$\bx\notin S_1(\bom)\implies \ [h_{\omega_1}\otimes t, \bx]v_\mu=0.$$ Set $$\bom_1=\sum_{\bx\in S_1(\bom)}\!\!\!\!z_\bx(x_{\beta_s}^-\otimes t^{\ell_s})\cdots (x_{\beta_1}^-\otimes t^{\ell_1+1})$$ and note that $\bom_1\in \bs$ is non--zero. Moreover,  $$[h_{\omega_1}\otimes t, \bom] v_\mu=\bom_1 v_\mu=0.$$ Since the minimal grade of $\bom_1$ is at least $k+1$,  the proposition is proved when $S_1(\bom)\neq \emptyset$.

Suppose now that $S_1(\bom)=\emptyset$, i.e., $$z_\bx\ne 0\implies \bx= (x_{\beta_s}^-\otimes t^{\ell_s})\cdots (x_{\beta_2}^-\otimes t^{\ell _2} )(x_{\beta_1}^-\otimes t^{N_{\beta_1}}).$$ If $\bom= x^-_{i_p.i_{p+1}}\otimes t^{\nu(h_{i_p,i_{p+1}})+1}),$ there is nothing to prove.
Otherwise, there exists  $i_p< j\le i_{p+1}$ minimal so that the set $$S_2(\bom)=\{\bx= (x_{\beta_s}^-\otimes t^{\ell_s})\cdots (x_{\beta_2}^-\otimes t^{\ell _2})(x_{\beta_1}^-\otimes t^{N_{\beta_1}})\in S: z_\bx\ne 0,\ \ \beta_1=\alpha_{i_p,j-1}\}\ne\emptyset.$$
 Then, $$[h_{\omega_j}\otimes t, \bx]=  (x_{\beta_s}^-\otimes t^{\ell_s})\cdots (x_{\beta_2}^-\otimes t^{\ell _2+1})(x_{\beta_1}^-\otimes t^{N_{\beta_1}}),\ \ \bx\in S_2(\bom), $$
and, using \eqref{e:reptcom}, we have $$[h_{\omega_j}\otimes t, \bx]v_\mu=0,\ \ \bx\notin S_2(\bom), \ \ z_\bx\ne 0.$$

Let $$\bom_2=\sum_{\bx\in S_2(\bom)}z_\bx (x_{\beta_s}^-\otimes t^{\ell_s}) \cdots (x_{\beta_2}^-\otimes t^{\ell_2+1})(x_{i_p,j-1}^-\otimes t^{N_{\alpha_{i_p,j-1}}}).$$ The preceding discussion proves that $$[h_{\omega_j}\otimes t, \bom]v_\mu=\bom_2v_\mu=0,$$ and that $\bom_2$ has minimal grade at least $k+1$. However, it need not be true that $\bom_2\in S$; for instance if $\bx\in S_2(\bom)$ is such that  $\ell_2=N_{\beta_2}$. To address this issue, we define a further subset
$S_2'(\bom)$  of $S_2(\bom)$ consisting of elements $\bx$ with $\ell_2+1\le N_{\beta_2}$. Let  $\bom'\in\bs $ be defined by:
$$\bom'=\left(\sum_{\bx\in S'_2(\bom)}z_\bx(x_{\beta_s}^-\otimes t^{\ell_s})\cdots (x_{\beta_2}^-\otimes t^{\ell_2+1})\right)(x_{i_p, j-1}^-\otimes t^{N_{\alpha_{i_p,j-1}}}).$$ Setting $S_2''(\bom)= S_2(\bom)\setminus S_2(\bom')$, we   note that$$(\bom_2-\bom')v_\mu= \left(\sum_{\bx\in S''_2{(\bom)}}z_\bx (x_{\beta_s}^-\otimes t^{\ell_s})\cdots(x_{\beta_3}^-\otimes t^{\ell_3} )(x_{\beta_2+\beta_1}^-\otimes  t^{N_{\beta_1+\beta_2}})\right)v_\mu .$$ The expression in parentheses on the right hand side is an element of $\bs$ and we denote it by $\bom''$.
 Moreover, the elements $\bom'$ and $\bom''$ are clearly linearly independent elements of $\bs$ assuming that at least one of them is  non--zero.  To see that this is in fact the case, assume that $\bom'=0$, i.e., $S_2'(\bom)=\emptyset$ and $S_2''(\bom)=S_2(\bom)$. This means that the  elements   $$ (x_{\beta_s}^-\otimes t^{\ell_s})\cdots(x_{\beta_3}^-\otimes t^{\ell_3} )(x_{\beta_2+\beta_1}^-\otimes  t^{N_{\beta_1+\beta_2}}), \   \  \  \bx\in S_2(\bom),$$  are all distinct, thus linearly independent (recall that $\beta_1=\alpha_{i_p,j-1}$) and, hence, $\bom''\ne 0$. Therefore $$0\ne \tilde\bom:=\bom'+\bom''\in \bs, \ \ \tilde\bom v_\mu=0,$$ and the minimal graded component of $\tilde\bom$ is at least $k+1$. The proposition is proved.

\end{pf}

\end{document}